%% file: agt-5-52.tex
\documentclass{gtart_h}
\input agtout

\lognumber{52}
\volumenumber{5}
\volumeyear{2005}
\papernumber{52}
\pagenumbers{1325}{1364}
\received{14 January 2005} 
\accepted{20 September 2005}
\published{6 October 2005}
 
\usepackage{amsmath, amssymb, graphicx, color}

\def\figref#1{\hyperlink{#1anchor}{Figure~\ref*{#1}}}
\def\anchor#1{\noindent\hypertarget{#1anchor}{\smash{$\phantom{99}$}}\newline}

\newtheorem{theorem}{Theorem}[section] 
\newtheorem{proposition}[theorem]{Proposition} 
\newtheorem{lemma} [theorem]{Lemma} 
\newtheorem{corollary} [theorem] {Corollary}

\theoremstyle{definition}
\newtheorem{definition} [theorem]{Definition} 
\newtheorem{observation} [theorem]{Observation}

\newtheorem{remark}[theorem]{Remark}

\newtheorem{convention}[theorem]{Convention}
\newtheorem{terminology}[theorem] {Terminology}
\newtheorem{example} [theorem] {Example}

\def\F{\mathcal F}
\def\R{\mathbb R}
\def\Z{\mathbb Z}

\def\SK {\underrightarrow{\text{Ker}}\, }
\def\Gn {\mathcal G_n}

\def\Finf{\mathcal F_\omega}

\def\SegCons{6\phi(4\delta)}
\def\A{\overline{A}}
\def\B{\overline{B}}
\def\C{\overline{C}}

\begin{document}

\title{Limits of (certain) CAT(0) groups,\\I: Compactification}

\author[Daniel Groves]{Daniel Groves}
\address{Department of Mathematics, California Institute of 
Technology\\Pasadena, CA, 91125, USA} 
\email{groves@caltech.edu} 


\asciikeywords{CAT(0) spaces, isolated flats, Limit groups, R-trees}

\begin{abstract}
The purpose of this paper is to investigate torsion-free groups which act properly and cocompactly on CAT$(0)$ metric spaces which have isolated flats, as defined by Hruska \cite{Hruska}.  Our approach is to seek results analogous to those of Sela, Kharlampovich and Miasnikov for free groups and to those of Sela (and Rips and Sela) for torsion-free hyperbolic groups.

This paper is the first in a series.  In this paper we extract an $\R$-tree from an asymptotic cone of certain CAT$(0)$ spaces.  This is analogous to a construction of Paulin, and allows a great deal of algebraic information to be inferred, most of which is left to future work.
\end{abstract}

\asciiabstract{%
The purpose of this paper is to investigate torsion-free groups which
act properly and cocompactly on CAT(0) metric spaces which have
isolated flats, as defined by Hruska.  Our approach is to seek results
analogous to those of Sela, Kharlampovich and Miasnikov for free
groups and to those of Sela (and Rips and Sela) for torsion-free
hyperbolic groups.  This paper is the first in a series.  In this
paper we extract an R-tree from an asymptotic cone of certain CAT(0)
spaces.  This is analogous to a construction of Paulin, and allows a
great deal of algebraic information to be inferred, most of which is
left to future work.}

\maketitle

\section{Introduction}

Using the theory of isometric actions on $\R$-trees as a starting point, Sela has solved the isomorphism problem for hyperbolic groups (at least for torsion-free hyperbolic groups which do not admit a small essential action on an $\R$-tree \cite{SelaIso}, though he has a proof in the general torsion-free case), has proved that torsion-free hyperbolic groups are Hopfian \cite{SelaHopf}, and recently has classified those groups with the same elementary theory as a given torsion-free hyperbolic group \cite{Sela1, SelaVI, SelaHyp}.  Kharlampovich and Miasnikov have a similar, but more combinatorial, approach to this last problem for free groups; see \cite{KM} and references contained therein.\footnote{Neither Sela's nor Kharlampovich and Miasnikov's work on the elementary theory of groups have entirely appeared in refereed journals.}

\eject

It seems that Sela's methods will not work for non-positively curved groups in general (whatever the phrase `non-positively curved group' means).  For example, Wise \cite{Wise} constructed a group which acts properly and cocompactly on a CAT$(0)$ metric space, but is non-Hopfian.

The class of groups acting properly and cocompactly on CAT$(0)$ spaces with the isolated flats condition is in many ways an intermediary between hyperbolic groups (which are the `negatively curved groups' in the context of discrete group) and CAT$(0)$ groups.  Sela \cite[Question I.8]{SelaProblems} asked whether such a group is Hopfian, and whether one can construct {\em Makanin-Razborov diagrams} for these groups.  In the second paper in this series, \cite{CWIFShort}, we will provide a positive answer to these questions (under certain extra hypotheses, described below).  The purpose of this paper is to develop tools for addressing these questions.

The initial ingredient in many of Sela's arguments is a result of Paulin (\cite{Paulin3, Paulin}; see also \cite{Bestvina} and \cite{BS}; and see \cite{MS84} for work preceding Paulin's) which extracts an isometric action on an $\R$-tree from (certain) sequences of actions on $\delta$-hyperbolic spaces.   Given two finitely generated groups $G$ and $\Gamma$, and a sequence of non-conjugate homomorphisms $\{ h_i \co G \to \Gamma \}$, it is straightforward to construct an action of $G$ on a certain asymptotic cone of $\Gamma$ with no global fixed point.  If $\Gamma$ acts properly and cocompactly by isometries on a metric space $X$, then a $G$-action can be constructed on an asymptotic cone of $X$ (which is bi-Lipschitz homeomorphic, but not necessarily isometric, to the analogous asymptotic cone of $G$).  For $\delta$-hyperbolic groups, this is in essence the above-mentioned result of Paulin.  In the case of groups acting on CAT$(0)$ spaces, it is carried out by Kapovich and Leeb in \cite{KL}, but the general case is hardly more complicated.  Of course, for a general finitely generated group $\Gamma$, the existence of an action of $G$ on an asymptotic cone of $\Gamma$ with no global fixed point provides little information about $G$ or $\Gamma$.  In this paper we place certain restrictions on $\Gamma$ so that we can find a $G$-action on an $\R$-tree which provides much the same information as Paulin's result.  We study the case where $\Gamma$ acts properly and cocompactly on a CAT$(0)$ metric space with isolated flats.

We study the asymptotic cone of such a space.  Under a further hypothesis (that the stabilisers of maximal flats are free abelian), we construct an $\R$-tree, which allows many of Sela's arguments to be carried out in this context (though we leave most such applications to subsequent work).

The first application of this construction is the following:\eject

\medskip

{\bf Theorem \ref{SplittingTheorem}}\qua
{\sl Suppose that $\Gamma$ is a torsion-free group acting properly
and cocompactly on a CAT$(0)$ space $X$ which has isolated flats,
so that flat stabilisers in $\Gamma$ are abelian.  Suppose further
that $\text{Out}(\Gamma)$ is infinite.  Then $\Gamma$ admits a nontrivial splitting over a finitely generated free abelian group.}

\medskip

This partially answers a question of Swarup (see \cite[Q 2.1]{BestvinaQuestions}).  However, Theorem \ref{SplittingTheorem} is only the first application.  Our hope is that much of Sela's program for free groups and torsion-free hyperbolic groups can be carried out for groups
$\Gamma$ as in the statement of Theorem \ref{SplittingTheorem}.  In future work, we will consider the automorphism groups of such groups (in analogy with \cite{RipsSelaGAFA, SelaGAFA}), the Hopf property (in analogy with \cite{SelaHopf}) and Makanin-Razborov diagrams for these groups (in analogy with \cite{Sela1, SelaHyp}).  The last of these involves finding a description of $\text{Hom}(G,\Gamma)$, where $G$ is an arbitrary finitely-generated group.   A key argument in Sela's solution to all of these problems for torsion-free hyperbolic groups is the {\em shortening argument}, which we present for these CAT$(0)$ groups
with isolated flats in \cite{CWIFShort}.

The outline of this paper is as follows.  In Section \ref{Prelim} we recall some basic definitions and results and prove some preliminary results about CAT$(0)$ spaces with isolated flats and groups acting properly and cocompactly on such spaces.  In Section \ref{Limits}, we consider a torsion-free group $\Gamma$ which acts properly and cocompactly on a CAT$(0)$ metric space $X$ with isolated flats. Given a finitely generated group $G$ and a sequence of homomorphisms $\{ h_n \co G \to \Gamma \}$ no two of which differ only by an inner automorphism of $\Gamma$, it is straightforward to construct an action of $G$ on the asymptotic cone of $X$.  A key feature of this action is that it has no global fixed point.    This construction amounts to a compactification of a certain space of $G$-actions on $X$ (those actions which factor through a fixed homomorphism $q \co \Gamma \to \text{Isom}(X)$).  In Section \ref{TreeSection}, we restrict to a torsion-free group $\Gamma$ which acts properly and cocompactly on a CAT$(0)$ space with isolated flats and has abelian flat stabilisers.
Under this additional hypothesis, we are able to extract an isometric action of $G$ on an $\R$-tree $T$ with no global fixed point.  The action of $G$ on $T$ largely encodes the same information from the homomorphisms $\{ h_n \}$ as Paulin's construction does in the case where $\Gamma$ is $\delta$-hyperbolic.  See, in particular, Theorem \ref{LinfProps}, the main technical result of this paper.  Finally, in Section \ref{conclusion} we discuss a few simple relations between our limiting objects, $\Gamma$-limit groups, and other definitions of $\Gamma$-limit groups, and prove Theorem \ref{SplittingTheorem}.  

I would like to thank Jason Manning for several conversations which illustrated my na\"ivet\'e, and in particular for pointing out an incorrect argument in a previous construction of the limiting tree $T$ in \S\ref{TreeSection}.

\section{CAT$(0)$ metric spaces with isolated flats and isometric actions 
upon them} \label{Prelim}

For the definition of $\R$-trees and the basic properties of their isometries, we refer the reader to \cite{Shalen1}, \cite{Shalen2}, \cite{Chiswell} and \cite{Bestvina2}.  For this paper, we do not need much of this theory.

For the definition and a multitude of results about CAT$(0)$ metric spaces, and isometric actions upon them,  we refer the reader to \cite{BH}.  We recall only a few basic properties and record our notation.

Suppose that $X$ is a geodesic metric space.  If $p,q,r \in X$, then $[p,q]$ denotes a geodesic between $p$ and $q$, and $\Delta (p,q,r)$ denotes the triangle consisting of the geodesics $[p,q],[q,r],[r,p]$.  Geodesics (and hence geodesic triangles) need not be unique in geodesic metric spaces, but they are in CAT$(0)$ spaces.  If $p,q,r \in X$ then $[p,q,r]$ denotes the path $[p,q] \cup [q,r]$.  Expressions such as $[p,q,r,s]$ are defined similarly.

If $\Gamma$ is a group acting properly and cocompactly by homeomorphisms on a connected simply-connected topological space then $\Gamma$ is finitely presented (see \cite[Theorem I.8.10, pp.135-137]{BH}). Obviously, if $\Gamma$ is torsion-free, then the action is free.

Suppose now that $X$ is a CAT$(0)$ metric space and that $\Gamma$ acts properly and cocompactly by isometries on $X$.  Then (see \cite[II.6.10.(2), p.233]{BH}) each element of $\Gamma$ acts either elliptically (fixing a point) or hyperbolically (there is an invariant axis upon which the element acts by translation).  If also $\Gamma$ is torsion-free then all isometries are hyperbolic.

Recall the following two results.

\begin{lemma} \label{Convex} {\rm\cite[Proposition II.2.2, p. 176]{BH}}\qua
Let $X$ be a CAT$(0)$ space.  Given any pair of geodesics $c\co [0,1] \to X$ and $c' \co [0,1] \to X$ parametrised proportional to arc length, the following inequality holds for all $t \in [0,1]$:
\[	d_X(c(t),c'(t)) \le (1-t)d_X(c(0),c'(0)) + t(d_X(c(1),c'(1)) .	\]
\end{lemma}

\begin{proposition} {\rm\cite[Proposition II.2.4, pp. 176--177]{BH}}\qua
Let $X$ be a CAT$(0)$ space, and let $C$ be a convex subset which is complete in the induced metric.  Then,
\begin{enumerate}
\item for every $x \in X$, there exists a unique point $\pi_C(x) \in C$ such that $d(x,\pi_C(x)) = d(x,C) := \inf_{y \in C}d(x,y)$;
\item if $x'$ belongs to the geodesic segment $[x,\pi_C(x)]$, then $\pi_C(x') = \pi_C(x)$;
\item the map $x \to \pi_C(x)$ is a retraction of $X$ onto $C$ which does not increase distances.
\end{enumerate}
\end{proposition}

\subsection{CAT$(0)$ spaces with isolated flats and groups acting on them}

\begin{definition}
A {\em flat} in a CAT$(0)$ space $X$ is an isometric embedding of Euclidean space $\mathbb E^k$ into $X$ for some $k \ge 2$.
\end{definition}
Note that we do not consider a geodesic line to be a flat.

\begin{definition} \cite[2.1.2]{Hruska} \label{CWIFDef}
A CAT$(0)$ metric space $X$ has {\em isolated flats} if it contains a family $\mathcal F_X$ of flats with the following properties:
\begin{enumerate}
\item (Maximal)\qua  There exists $B \ge 0$ such that every flat in $X$ is contained in a $B$-neighbourhood of some flat in $\mathcal F_X$;
\item (Isolated)\qua There is a function $\phi \co \R_+ \to \R_+$ such that for every pair of distinct flats $E_1, E_2 \in \F_X$ and for every $k \ge 0$, the intersection of the $k$-neighbourhoods of $E_1$ and $E_2$ has diameter less than $\phi(k)$.
\end{enumerate}
\end{definition}

This definition is due to C. Hruska \cite{Hruska}, but such an idea is implicit in Chapter 11 of \cite{E+}, and in the work of Wise \cite{Wise2} and of Kapovich and Leeb \cite{KL}.

\begin{convention} \label{PhiConvention}
To simplify constants in the sequel, we assume that $\phi(k) \ge k$ for all $k \ge 0$ and that $\phi$ is a nondecreasing function.  We can certainly make these assumptions, and usually do so without comment.
\end{convention}

For the basic properties of CAT$(0)$ metric spaces with isolated flats, for examples of such spaces, and for some properties of isometric actions upon them, we refer the reader to \cite{Hruska}.

\begin{figure}[ht!]\anchor{RelThinPic}
\begin{center}
\scalebox{0.9}{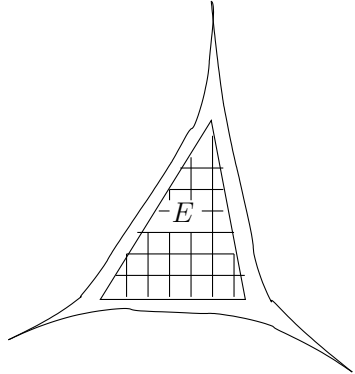}
\end{center}
\caption{A triangle which is thin relative to the flat $E$}
\label{RelThinPic}
\end{figure}

Hruska also introduced the {\em relatively thin triangles} property:

\begin{definition} \cite[3.1.1]{Hruska}\label{RelThinDef}\qua
A geodesic triangle in a metric space $X$ is {\em $\delta$-thin relative to the flat $E$} if each side of the triangle lies in the $\delta$-neighbourhood of the union of $E$ and the other two sides of the triangle (see \figref{RelThinPic}).  A metric space $X$ has the {\em relatively thin triangle property} if there is a constant $\delta$ so that each triangle in $X$ is either $\delta$-thin in the usual sense or $\delta$-thin relative to some flat in $\mathcal F_X$.
\end{definition}

Using work of Dru\c{t}u and Sapir \cite{DS} on asymptotic cones of relatively hyperbolic groups, Hruska and Kleiner \cite{HruskaKleiner} have proved that if $X$ is a CAT$(0)$ space with isolated flats which admits a cocompact isometric group action then $X$ satisfies the relatively thin triangles condition.  In this paper,
the symbol `$\delta$' will always refer to the constant from Definition
\ref{RelThinDef}.

\begin{terminology}
When we refer to a {\em CAT$(0)$ group with isolated flats} we mean
a group which admits a proper, cocompact and isometric action on a CAT$(0)$ space with isolated flats.
\end{terminology}

We now consider some of the basic properties of CAT$(0)$ spaces with isolated flats, and groups acting properly, cocompactly and isometrically upon them, which are necessary in the sequel.

\begin{proposition} {\rm\cite[2.1.4]{Hruska}}\label{FInv}\qua
Suppose $X$ is a CAT$(0)$ space with isolated flats.  The family $\F_X$ of flats in Definition \ref{CWIFDef} may be assumed to be invariant under all isometries of $X$.
\end{proposition}

\begin{lemma} {\rm\cite[2.1.9]{Hruska}}\label{Periodic}\qua
Suppose that the CAT$(0)$ space $X$ has isolated flats and admits a proper and cocompact action by some group of isometries.  Then any maximal flat in $X$ is periodic.
\end{lemma}

\begin{lemma} \label{UniqueFlat}
Suppose that $X$ is a CAT$(0)$ space with isolated flats, and that $\Delta = \Delta(a,b,c)$ is a geodesic triangle in $X$.  If $\Delta$ is not $(\delta + \frac{\phi(\delta)}{2})$-thin then $\Delta$ is $\delta$-thin relative to a {\em unique} flat $E \in \F_X$.
\end{lemma}
\begin{proof}
Let $l_{a,b}$ be that part of the geodesic $[a,b]$ which lies outside of the $\delta$-neighbourhood of $[a,c] \cup [b,c]$, and define $l_{a,c}$ and $l_{b,c}$ similarly.  

Suppose that $\Delta$ is $\delta$-thin relative to $E, E' \in \F_X$, where $E \neq E'$.  Then $l_{a,b}, l_{a,c}, l_{b,c}$ all lie in the $\delta$-neighbourhood both of $E$ and of $E'$.  The intersection of these $\delta$-neighbourhoods has diameter at most $\phi(\delta)$.  Therefore, the length of $l_{a,b}$ is at most $\phi(\delta)$ (since it is a geodesic).  Thus, from any point on $l_{a,b}$, the distance to $[a,c] \cup [b,c]$ is at most $\delta + \frac{\phi(\delta)}{2}$.

A symmetric argument for $l_{a,c}$ and $l_{b,c}$ finishes the proof.
\end{proof}

\subsection{Bieberbach groups and toral actions on CAT$(0)$ spaces with isolated flats}

Given a proper and cocompact isometric action of a group $\Gamma$ on a CAT$(0)$ space $X$ with isolated flats, we are compelled to study the subgroups $\text{Stab}(E)$, where $E$ is a maximal flat in $X$.\footnote{By $\text{Stab}(E)$ we mean $\{ g \in \Gamma\ |\ g.E = E \}$.  The point-wise stabiliser is $\text{Fix}(E) = \{ g \in \Gamma\ |\ g.x = x, \ \forall x \in E \}$.}

By Lemma \ref{Periodic} we have a proper and cocompact action of the group $\text{Stab}_\Gamma(E)$ on $E \cong \mathbb E^n$.  Recall the following celebrated result of Bieberbach \cite{Bieb1, Bieb2}. \footnote{An {\em $n$-dimensional crystallographic group} is a cocompact discrete group of isometries of $\mathbb E^n$.}

\begin{theorem} [Bieberbach; see for example \cite{Thurston}, 4.2.2, p.222]
$\phantom{99}$
\begin{enumerate}
\item[\rm(a)] A group $\Gamma$ is isomorphic to a discrete group of isometries of $\mathbb E^n$, for some $n$, if and only if $\Gamma$ contains a subgroup of finite index that is free abelian of finite rank;
\item[\rm(b)] An $n$-dimensional crystallographic group $\Gamma$ contains a normal subgroup of finite index that is free abelian of rank $n$ and equals its own centraliser.  This subgroup is characterised as the unique maximal abelian subgroup of finite index in $\Gamma$, or as the translation subgroup of $\Gamma$.
\end{enumerate}
\end{theorem}

The structure of the subgroups $\text{Stab}_\Gamma(E)$ will be important to us in the sequel.  In particular, when there are such groups which are not free abelian the construction in Section \ref{TreeSection} does not work.  Motivated by this consideration, we make the following

\begin{definition}
Suppose that $X$ is a CAT$(0)$ space with isolated flats and that a  group $\Gamma$ acts properly and cocompactly by isometries on $X$.  We say that the action of $\Gamma$ on $X$ is {\em toral} if for each maximal flat $E \subseteq X$, the subgroup $\text{Stab}(E) \le \Gamma$ is free abelian.  We say that $\Gamma$ is a {\em toral} CAT$(0)$ group with isolated flats if there is a proper, cocompact and toral action of $\Gamma$ on a CAT$(0)$ space $X$ with isolated flats.
\end{definition}

\begin{remark}
We observe in Lemma \ref{AllToral} below that if a torsion-free group $\Gamma$ admits a proper, cocompact and toral action on some CAT$(0)$ space $X$ with isolated flats then {\em any} proper and cocompact action of $\Gamma$ on a CAT$(0)$ space with isolated flats is toral.  Thus the property of being toral belongs to the group rather than the given action on a CAT$(0)$ space with isolated flats.  Also, Hruska and Kleiner have proved \cite{HruskaKleiner} that any CAT$(0)$ space $X$ on which a CAT$(0)$ group with isolated flats acts properly and cocompactly by isometries has isolated flats.
\end{remark}

\subsection{Basic algebraic properties of CAT$(0)$ groups with isolated flats}

In this paragraph we consider a few basic algebraic properties of torsion-free CAT$(0)$ groups with isolated flats.

\begin{definition}
A subgroup $K$ of a group $G$ is said to be {\em malnormal} if for all $g \in G \smallsetminus K$ we have $gKg^{-1} \cap K = \{ 1 \}$.  

A group $G$ is said to be {\em CSA} if any maximal abelian subgroup of $G$ is malnormal.
\end{definition}

The following lemma is straightforward and certainly well known, but we record and prove it for later use.

\begin{lemma} \label{SolAb}
Suppose that $G$ is a CSA group.  Then every soluble subgroup of $G$ is abelian.  Also, every virtually abelian subgroup of $G$ is abelian.
\end{lemma}
\begin{proof}
Suppose that $S$ is a nontrivial soluble subgroup of $G$.  Let $S^{(i)}$ be the smallest nontrivial term of the derived series of $S$.  Then $S^{(i)}$ is a normal abelian subgroup of $S$.  However, it is an abelian subgroup of $G$, so is contained in a maximal abelian subgroup $A$.  If $g \in S$, then $g$ normalises $S^{(i)}$, so $g \in A$, since $A$ is malnormal.  Therefore, $S$ is contained in $A$ and $S$ is abelian.

Any virtually abelian subgroup $H$ has a finite index normal abelian subgroup $A$.  By the above argument, the normaliser of $A$ is abelian and contains $H$, so $H$ is abelian.
\end{proof}

\begin{proposition} \label{malnormal}
Suppose that $\Gamma$ is a torsion-free group which admits a proper and cocompact action on a CAT$(0)$ space $X$ with isolated flats.  Then the stabiliser in $\Gamma$ of any maximal flat in $X$ is malnormal.
\end{proposition}
\begin{proof}
Let $\mathcal F_X$ be the collection of flats from Definition \ref{CWIFDef}, and let $E$ be a maximal flat in $X$.  Consider $M = \text{Stab}(E)$.  Without loss of generality, we may assume that $E \in \mathcal F_X$.

Suppose that $g \in \Gamma$ is such that $gMg^{-1} \cap M \neq \{ 1 \}$.  We prove that $g \in M$.  There exist $a_1 , a_2 \in M \smallsetminus \{ 1 \}$ so that $ga_1g^{-1} = a_2$.

Now, $ga_1g^{-1} = a_2$ leaves  both $E$ and $gE$ invariant.  Therefore, there is an axis for $ga_1g^{-1}$ in each of $E$ and $gE$, and there is a Euclidean strip, isometric to $[0,k] \times \R$ for some $k$, joining these axes.  However, $E$ and $gE$ are both in $\mathcal F_X$ by Proposition \ref{FInv}, and we have seen that the $k$-neighbourhoods of $E$ and $gE$ intersect in an unbounded set, so we must have that $E = gE$, which is to say that $g \in M$.
\end{proof}

\begin{corollary} \label{CSA}
Suppose that $\Gamma$ is a torsion-free toral CAT$(0)$ group with isolated flats.  Then $\Gamma$ is CSA.
\end{corollary}
\begin{proof}
Let $A$ be a maximal abelian subgroup of $\Gamma$ and let $X$ be a CAT$(0)$ space with isolated flats with a proper, cocompact and toral action of $\Gamma$.

Suppose first that $A$ is noncyclic.  Then $A$ stabilises some flat $E \in \F_X$, and hence some maximal flat (by the Isolated Flats condition).  Since $A$ is maximal abelian, and the action of $\Gamma$ on $X$ is toral, $A = \text{Stab}(E)$.  In this case the result follows from Proposition \ref{malnormal}.

Suppose now that $A$ is a cyclic maximal abelian subgroup, and that for some $g \in \Gamma \smallsetminus A$ we have $gAg^{-1} \cap A \neq \{ 1 \}$.  Let $A = \langle a \rangle$.  Then $ga^pg^{-1} = a^q$ for some $p, q$.  Since $A$ is maximal abelian, we do not have $p,q = 1$.  However, $\Gamma$ is a CAT$(0)$ group, so $|p| = |q|$ (see \cite[Theorem III.$\Gamma$.1.1(iii)]{BH}).  Thus $g^2$ commutes with $a^p$.  Therefore, $\langle a^p \rangle$ is central in $G = \langle g^2, a^p \rangle$.  By \cite[II.6.12]{BH}, there is a finite index subgroup $H$ of $G$ so that $H = \langle a^p \rangle \times H_1$, for some group $H_1$.  If $H_1$ is infinite, then $\langle a^p \rangle$ is contained in a subgroup isomorphic to $\mathbb Z^2$.  This $\mathbb Z^2$ stabilises a flat, and hence a maximal flat, so $\langle a^p \rangle$ is contained in $\text{Stab}(E)$ for some $E \in \F_X$.  However, $a$ normalises $\langle a^p \rangle$, so by Proposition \ref{malnormal} $a \in \text{Stab}(E)$.  This subgroup is abelian since the action of $\Gamma$ on $X$ is toral, which contradicts $A$ being maximal abelian.  Therefore, $H_1$ is finite, and since $\Gamma$ is torsion-free, $H_1$ is trivial.  Therefore, $G$ is virtually cyclic, and being torsion-free, is itself infinite cyclic.  Hence $g^2$ commutes with $a$ and so $\langle g^2 \rangle$ is central in $G_1 = \langle g, a \rangle$.  Exactly the same argument as above applied to $G_1$ and $\langle g^2 \rangle$ implies that $G_1$ is cyclic.  Since $A$ is maximal abelian, $g \in A$, a contradiction to the choice of
$g$.  Thus $A$ is malnormal, as required.
\end{proof}

\begin{lemma} \label{AllToral}
Suppose that $\Gamma$ is a torsion-free group which admits a proper, cocompact and toral action on a CAT$(0)$ space $X$ with isolated flats.  Then any proper and cocompact action of $\Gamma$ on a CAT$(0)$ space with isolated flats is toral.  If $\Gamma$ is a torsion-free CAT$(0)$ group with isolated flats then $\Gamma$ is toral if and only if $\Gamma$ is CSA.
\end{lemma}
\begin{proof}
Let $\Gamma$ act properly and cocompactly on a CAT$(0)$ space $Y$ with isolated flats, and let $M$ be the stabiliser of a maximal flat $E \in \F_Y$.

Since $\Gamma$ admits a proper, cocompact and toral action on a CAT$(0)$ space $X$ with isolated flats, by Corollary \ref{CSA} any maximal abelian subgroup of $\Gamma$ is malnormal, and so the normaliser of any abelian group is abelian.

Since $M$ is a Bieberbach group, it has a normal abelian subgroup $A$ of finite index.  However, by the above, the normaliser of $A$ is abelian and it certainly contains $M$, so $M$ is abelian.  Therefore the action of $\Gamma$ on $Y$ is toral.  This proves the first claim of the lemma.  The second claim follows from the proof of the first and Corollary \ref{CSA}.
\end{proof}

\begin{definition}
A group $G$ is said to be {\em commutative transitive} if for all $u_1,u_2,u_3 \in G \smallsetminus \{ 1 \}$, whenever $[u_1,u_2] = 1$ and $[u_2,u_3] = 1$ we necessarily have $[u_1,u_3] =1$.
\end{definition}

CSA groups are certainly commutative transitive, so we have

\begin{corollary} \label{CommTrans}
Suppose $\Gamma$ is torsion-free toral CAT$(0)$ group with isolated flats.  Then $\Gamma$ is commutative transitive.  Hence every abelian subgroup in $\Gamma$ is contained in a unique maximal abelian subgroup.
\end{corollary}

\subsection{Projecting to flats} \label{ProjectSection}

Fix $X$, a CAT$(0)$ space with isolated flats.  Let $\delta$ be the constant from Definition \ref{RelThinDef} and $\phi$ the function from Definition \ref{CWIFDef}.

We study the closest-point projection from $X$ onto a flat $E \subset X$.

\begin{lemma} \label{Triangle}
Suppose that $E \in \mathcal F_X$.  Suppose that $x,y \in E$ and $z \in X$.  There exist $u \in [x,z]$ and $v \in [y,z]$ so that $u$ and $v$ lie in the $2\delta$-neighbourhood of $E$, and 
\[	d_X(u,v) \le \phi(\delta) .	\]
\end{lemma}
\begin{proof}
If $z$ lies in the $2\delta$-neighbourhood of $E$ then the result is immediate, so we assume that this is not the case.  The key (though trivial) observation is that $[x,y]$ lies entirely within $E$.

Let $u_1$ be the point on $[x,z]$ which is furthest from $x$ in the $\delta$-neighbourhood of $E$.  The convexity of $E$ and the convexity of the metric on $X$ ensures that $u_1$ is unique.

We consider the triangle $\Delta = \Delta(x,y,z)$.  If $\Delta$ is $\delta$-thin, then there is clearly a point $v_1$ on $[y,z]$ within $\delta$ of $v_1$, and we take $u = u_1, v = v_1$ (this is because if $u_1$ is not $\delta$-close to $[y,z]$ then neither is a point nearby $u_1$, but there are points arbitrarily close to $u_1$ on $[x,z]$ which are not $\delta$-close to $E$, which would contradict $\Delta$ being $\delta$-thin since $[x,y] \subset E$).

Thus suppose that $\Delta$ is not $\delta$-thin, so that $\Delta$ is $\delta$-thin relative to a flat $E'$.  If $E' = E$, then we have a point $v_1 \in [y,z]$ which is within $\delta$ of $u_1$, by the same reasoning as above.  Again, we take $u = u_1, v = v_1$.  

Suppose then that $E' \ne E$.  Either $u_1$ is $\delta$-close to $[y,z]$, in which case we proceed as above, or $u_1$ is $\delta$-close to $E'$.  In this case define $v_2$ to be the point on $[y,z]$ which is furthest from $y$ but in the $\delta$-neighbourhood of $E$.  Again, $v_2$ is either $\delta$-close to $[x,z]$ or $\delta$-close to $E'$.  In the first situation, we proceed as above, with $v = v_2$ and $u$ a point on $[x,z]$ which is within $\delta$ of $v_2$.  In the second situation, both $u_1$ and $v_2$ are within $\delta$ of $E$ and of $E'$, and the intersection of the $\delta$-neighbourhoods of $E$ and $E'$ has diameter less than $\phi(\delta)$.  Thus in this case, $d_X(u_1,v_2) < \phi(\delta)$ and we may take $u = u_1, v = v_2$.

We have proved that there exist $u$ and $v$ in the $2\delta$-neighbourhood of $E$ so that $d_X(u,v) \le \max \{ \delta , \phi(\delta) \}$.  However, $\delta \le \phi(\delta)$ by Convention \ref{PhiConvention}, so the proof is complete.
\end{proof}

\begin{proposition} \label{Projection}
Suppose that $E \subset \mathcal F_X$ is a flat and that $x,y \in X$ are such that $[x,y]$ does not intersect the $4\delta$-neighbourhood of $E$.  Let $\pi \co X \to E$ be the closest-point projection map.  Then $d_X(\pi(x),\pi(y)) \le 2\phi(3\delta)$.
\end{proposition}
\begin{proof}
By Lemma \ref{Triangle}, there are points $w_1 \in [\pi(x),y]$ and $w_2 \in [\pi(y),y]$, both in the $2\delta$-neighbourhood of $E$, so that $d_X(w_1,w_2) \le \phi(\delta)$.

Now consider the triangle $\Delta' = \Delta(\pi(x),x,y)$.  By a similar argument to the proof of Lemma \ref{Triangle}, we find points $u_1 \in [\pi(x),x]$ and $u_2 \in [\pi(x),y]$ which lie outside the $2\delta$-neighbourhood of $E$ such that $d_X(u_1,u_2) \le \max \{ \delta, \phi(3\delta) \} \le \phi(3\delta)$.  Indeed, let $v_1$ be the point on $[\pi(x),x]$ furthest from $\pi(x)$ which lies in the $3\delta$-neighbourhood of $E$, and let $v_2$ be the point on $[\pi(x),y]$ furthest from $\pi(y)$ which lies in the $3\delta$-neighbourhood of $E$.  

If $\Delta'$ is $\delta$-thin then there is a point $u_2$ on $[\pi(x),y]$ within $\delta$ of $v_1$.  We may take $v_1 = u_1$.  Similarly, if $\Delta'$ is $\delta$-thin relative to $E$ then once again there must be such a point $u_2$.

Therefore, suppose that $\Delta'$ is $\delta$-thin relative to $E' \ne E$.  Then $v_1$ does not lie within $\delta$ of $[x,y]$ since $[x,y]$ does not intersect the $4\delta$-neighbourhood of $E$.  Therefore, either $v_1$ lies within $\delta$ of $E'$ or within $\delta$ of $[\pi(y),x]$.  The second case is unproblematic as usual.  Also, $v_2$ either lies within $\delta$ of $E'$ or within $\delta$ of $[\pi(x),x]$, and in this second case we proceed as usual.

So suppose that $v_1$ and $v_2$ both lie within $\delta$ of $E'$.  Then they both lie within the $3\delta$-neighbourhoods of $E$ and $E'$ and so $d_X(v_1,v_2) \le \phi(3\delta)$.

Now, $u_2$ is closer to $y$ along $[\pi(x),y]$ than $w_1$, since $u_2$ lies outside the $2\delta$-neighbourhood of $E$, and $w_1$ lies within.  Hence, the convexity of the metric in $X$ ensures that there is a point $u_3 \in [\pi(y),y]$ so that $d_X(u_2,u_3) \le d_X(w_1,w_2)$.

Now, $u_1 \in [\pi(x),x]$ so $\pi(u_1) = \pi(x)$, and similarly $\pi(y_3) = \pi(y)$.  Therefore,
\begin{eqnarray*}
d_X(\pi(x),\pi(y)) & = & d_X(\pi(u_1),\pi(u_3))\\
& \le & d_X(\pi(u_1),\pi(u_2)) + d_X(\pi(u_2),\pi(u_3))\\
& \le & d_X(\pi(u_1),\pi(u_2)) + d_X(\pi(w_1),\pi(w_2))\\
& \le & \phi(3\delta) + \phi(\delta) \\
& \le & 2\phi(3\delta),
\end{eqnarray*}
as required.
\end{proof}

\section{Asymptotic cones of CAT$(0)$ spaces with isolated flats} \label{Limits}

In this section, we construct a limiting action from a sequence of homomorphisms from a fixed finitely generated group $G$ to $\Gamma$, a CAT$(0)$ group with isolated flats.  The action is of $G$ on the asymptotic cone of $X$, where $X$ is the CAT$(0)$ space with isolated flats upon which $G$ acts properly and cocompactly.  Asymptotic cones of CAT$(0)$ spaces have been studied in \cite{KL} and we use or adapt many of their results.  We note that one of the results of this section is that the asymptotic cone of $X$ is a {\em tree-graded metric space}, in the terminology of \cite{DS}.  This follows from \cite{HruskaKleiner} and \cite{DS}.  This paper was written before \cite{HruskaKleiner} or \cite{DS} appeared publicly, and we need more results from this section than follow directly from either \cite{DS} or \cite{HruskaKleiner}.  Thus, we
prefer to leave this section unchanged, rather than referring to \cite{DS}
or \cite{HruskaKleiner} for some of the results herein. 

\begin{remark}
The construction below could be carried out in a similar way to those found in \cite{Paulin3, Paulin} (see also \cite{Bestvina} and \cite{BS}) using the equivariant Gromov topology on `approximate convex hulls' of finite orbits of a basepoint $x$ under the various actions of $G$ on $X$.  For the sake of brevity, however, we use asymptotic cones.  However, having used asymptotic cones we use Lemma \ref{GromovTop} below to pass back to the context of the equivariant Gromov topology.
\end{remark}

\subsection{Constructing the asymptotic cone}

Suppose that $X$ is a CAT$(0)$ space with isolated flats, and $\Gamma \to \text{Isom}(X)$ is a proper, cocompact and isometric action of $\Gamma$ on $X$.

Let $G$ be a finitely generated group, and suppose that  $\{ h_n \co G \to \Gamma \}$ is a sequence of nontrivial homomorphisms.  A homomorphism $h \co G \to \Gamma$ gives rise to a sequence of proper isometric actions of $G$ on $X$:
\[	\lambda_{h} \co G \times X,	\]
given by $\lambda_{h} = \iota \circ h$, where $\iota\co \Gamma \to \text{Isom}(X)$ is the fixed homomorphism given by the action of $\Gamma$ on $X$.

Because the action of $\Gamma$ on $X$ is proper and cocompact, 
we have the following:

\begin{lemma}
For any $y \in X$, $j \ge 1$ and $g \in G$, the function $\iota_{g,j,y} \co \Gamma \to \R$ defined by
\[	\iota_{g,j,y}(\gamma) = d_X \left( \gamma . y, \lambda_j ( g , \gamma . y ) \right) ,	\]
achieves its infimum for some $\gamma' \in \Gamma$.
\end{lemma}

Let $\mathcal A$ be a finite generating set for $G$ and let $x \in X$ be arbitrary.  For a homomorphism $h \co G \to \Gamma$ define $\mu_h$ and $\gamma_h \in \Gamma$ so that
\begin{eqnarray*}
\mu_h & = & \max_{g \in \mathcal A} d_X \left( \gamma_h . x , \lambda_h (g, \gamma_h . x) \right) \\
& = & \min_{\gamma \in \Gamma} \max_{g \in \mathcal A} d_X \left( \gamma . x , \lambda_h (g, \gamma . x) \right) .
\end{eqnarray*}
For the chosen sequence of homomorphisms $h_n \co G \to \Gamma$, we write $\lambda_n$ instead of $\lambda_{h_n}$, $\mu_i$ instead of $\mu_{h_i}$ and $\gamma_i$ instead of $\gamma_{h_i}$.

Now define the pointed metric spaces $(X_n,x_n)$ to be the set $X$ with basepoint $x_n = x$, with the metric $d_{X_n} = \frac{1}{\mu_n}d_X$.  Since there is a natural identification between $\text{Isom}(X)$ and $\text{Isom}(X_i)$, we consider $\lambda_i$ to give an action of $G$ on $X_i$, as well as on $X$.

The next lemma follows from the fact that $\Gamma . x$ is discrete, and that $G$ is finitely generated.

\begin{lemma} \label{mujGrows}
Suppose that for all $j \ne i$ there is no element $\gamma \in \Gamma$ so that $h_i = \tau_{\gamma} \circ h_j$ where $\tau_\gamma$ is the inner automorphism of $\Gamma$ induced by $\gamma$.  Then the sequence $\{ \mu_j \}$ does not contain a bounded subsequence.
\end{lemma}

We use the homomorphism $h_j$ and the translation minimising element $\gamma_j$ to define an isometric action $\hat{\lambda}_j \co G \times X_n \to X_n$ by defining
\[	\hat{\lambda}_n(g,y) = \left( \gamma_n^{-1} h_n(g) \gamma_n \right) . y	.	\]

\begin{convention}
For the remainder of the paper, we assume that the homomorphisms $h_n$ were chosen so that $\gamma_n = 1$ for all $n$.  Therefore, $\hat{\lambda}_n (g,x) = \lambda_n(g,x) = h_n(g) . x$ for all $n \ge 1$, $g \in G$ and $x \in X$.
\end{convention}

Using the spaces $(X_n,x_n)$ and the actions $\lambda_n$ of $G$ on $X_n$, we construct an action of $G$ on the asymptotic cone of $X$, with respect to the basepoints $x_n = x$, scalars $\mu_n$ and an arbitrary non-principal ultrafilter $\omega$.

We briefly recall the definition of asymptotic cones.  For more details, see \cite{VW} and \cite{DS}, or \cite{KL} in the context of  CAT$(0)$ spaces.

\begin{definition}
A {\em non-principal ultra-filter}, $\omega$, is a $\{ 0, 1 \}$-valued finitely additive measure on $\mathbb N$ defined on all subsets of $\mathbb N$ so that any finite set has measure $0$.
\end{definition}

The existence of non-principal ultrafilters is guaranteed by Zorn's Lemma.

Fix once and for all a non-principal ultrafilter $\omega$.\footnote{The choice of ultrafilter will affect the resulting construction, but will not affect our results.  Thus we are unconcerned which ultrafilter is chosen.}  Given any bounded sequence $\{ a_n \} \subset \mathbb \R$ there is a unique number $a \in \R$ so that, for all $\epsilon > 0$, $\omega \left( \{ a_n\ |\ |a - a_n| < \epsilon \} \right) = 1$.  We denote $a$ by $\omega$-$\lim \{ a_n \}$.  This notion of limit exhibits most of the properties of the usual limit (see \cite{VW}).

The {\em asymptotic cone of $X$ with respect to $\{ x_n \}$, $\{ \mu_n \}$ and $\omega$}, denoted $X_\omega$ is defined as follows.  First, define the set $\tilde{X}_\omega$ to consist of all sequences $\{ y_n \ | \ y_n \in X_n\}$ for which $\{ d_{X_n}(x_n,y_n) \}$ is a bounded sequence.  Define a pseudo-metric $\tilde{d}$ on $\tilde{X}$ by
\[	\tilde{d} (\{y_n, z_n \}) = \mbox{$\omega$-$\lim$ } d_{X_n}(y_n,z_n).	\]

The asymptotic cone $X_\omega$ is defined to be the metric space induced by the pseudo-metric $\tilde{d}$ on $\tilde{X}_\omega$:
\[	X_\omega := \tilde{X_\omega}/ \sim ,	\]
where the equivalence relation `$\sim$' on $\tilde{X}_\omega$ is defined by: $x \sim y$ if and only if $\tilde{d}(x,y) = 0$.  The pseudo-metric $\tilde{d}$ on $\tilde{X}_\omega$ naturally descends to a metric $d_\omega$ on $X_\omega$.

\begin{lemma} [see \cite{VW} and \cite{KL}, Proposition 3.6]
$(X_\omega, d_\omega)$ is a complete, geodesic CAT$(0)$ space.
\end{lemma}

We now define an isometric action of $G$ on $X_\omega$.  Let $g \in G$ and $\{ y_n \} \in \tilde{X}_\omega$.  Then define $g.\{ y_n \}$ to be $\{ {\lambda}_n(g, y_n) \} \in \tilde{X}_\omega$.  This descends to an isometric action of $G$ on $X_\omega$. 

\begin{remark}
The action of $G$ defined on the asymptotic cone $X_\omega$ is slightly different to the one described in \cite[\S 3.4]{KL}, but the salient features remain the same.\footnote{The difference comes in the choice of scalars $\mu_n$ and the choice of basepoints $x_n$.}
\end{remark}
 
We assume a familiarity with asymptotic cones, but we essentially use only two properties.  The first is that finite sets in $\mathbb N$ have $\omega$-measure $0$.  The second property is the following:

\begin{lemma}
Suppose that $X_\omega$ is constructed using the sequence $\{ h_n \co G \to \Gamma \}$ as above.  Suppose that $Q \subset G$ is finite and $S \subset X_\omega$ is finite.  For each $s \in S$, let $\{ s_n \}$ be a sequence of elements from $X_n$ such that $\{ s_n \} \in \tilde{X}_\omega$ is a representative of the equivalence class $s$.  Fix $\epsilon > 0$ and define $I_{\epsilon, Q, S}$ to be the set of $i \in \mathbb N$ so that for all $q_1, q_2 \in Q \cup \{ 1 \}$ and all $s, s' \in S$ we have 
\[	|d_{X_i}(\hat{\lambda}(q_1 , s_i),\hat{\lambda}(q_2 , s_i')) - d_{X_\omega}(q_1 . s, q_2 . s')| < \epsilon.  	\]
Then $\omega (I_{\epsilon,Q,S}) = 1$.
\end{lemma}

Given finite subsets $Q$ of $G$ and $S$ of $X_\omega$ and $\epsilon > 0$ as above, if $i \in I_{\epsilon,Q,S}$ then the pair $(X_i,\lambda_i)$ is called an {\em $\epsilon$-approximation} for $Q$ and $S$.

\subsection{Properties of $X_\omega$}

\begin{lemma} \label{NoXinfFixedPt}
The action of $G$ on $X_{\omega}$ by isometries does not have a global fixed point.
\end{lemma}
\begin{proof}
Let $K \subseteq X$ be a compact set so that the basepoint $x$ is in $K$ and $\Gamma . K = X$.  Let $D = \text{Diam}(K)$.  

Suppose that $y \in X_\omega$ is fixed by all points of $G$.  Choose a large $i$ so that (i) $\mu_i > 4D$; and (ii) $X_i$ is a $\frac{1}{2}$-approximation for $\{ y \}$ and $\mathcal A$.  (Recall that $\mathcal A$ is the fixed finite generating set for $G$.)  Thus, if $\{ y_n \}$ represents $y$ then for all $g \in \mathcal A$
\[	d_{X_i}(y_i,{\lambda}_i(g,  y_i)) < \frac{1}{2}.	\]
This implies that, for all $g \in \mathcal A$,
\[	d_X(y_i, h_i(g) . y_i) < \frac{\mu_i}{2}. 	\]
Now, there exists $\gamma \in \Gamma$ so that $d_X(x,\gamma . y_i) \le D$.  Let $g \in S$ be the element which realises the maximum:
\[ \mu_i  =  \min_{\overline{\gamma} \in \Gamma} \max_{g \in S} d_X(x, \overline{\gamma} h_i(g) \overline{\gamma}^{-1}) . x ).	\]
Then we have
\begin{eqnarray*}
\mu_i & \le & d_X(x,(\gamma h_i(g) \gamma^{-1}) . x) \\
& \le & d_X(x,\gamma . y_i) + d_X( \gamma . y_i, (\gamma h_i(g) \gamma^{-1}) \gamma . y_i) \\
&& + d_X((\gamma h_i(g) \gamma^{-1}) \gamma . y_i, (\gamma h_i(g) \gamma^{-1}) . x) \\
& = & 2d_X(x,\gamma . y_i) + d_X( y_i, h_i(g) . y_i)\\
& < & 2D + \frac{\mu_i}{2}.
\end{eqnarray*}
Since $\mu_i > 4D$ this is a contradiction.  Therefore there is no global fixed point for the action of $G$ on $X_\infty$.
\end{proof}

We now prove some results about $X_\omega$ which are very similar to those obtained in  \cite{KL} in the context of asymptotic cones of certain $3$-manifolds.

\begin{definition} (See \cite{KL}, $\S$2-2)\qua
Let $X$ be a CAT$(0)$ space and $x,y,z \in X$.  Define $x',y',z'$ by $[x,x'] = [x,y] \cap [x,z]$, $[y,y'] = [y,z] \cap [y,x]$ and $[z,z']  = [z,x] \cap [z,y]$.  The triangle $\Delta (x',y',z')$ is called the {\em open triangle} spanned by $x,y,z$.  The triangle $\Delta (x,y,z)$ is called {\em open} if $x=x'$, $y=y'$ and $z= z'$.  An open triangle $\Delta (x',y',z')$ is {\em non-degenerate} if the three points $x', y', z'$ are distinct.
\end{definition}

Let $\mathcal F_X$ be the set of flats in $X$ from Definition \ref{CWIFDef}.  Let $\mathcal F_n$ be the set $\mathcal F_X$ considered as subsets of $X_n$.  Denote by $\mathcal F_\omega$ the set of all flats in $X_\omega$ which arise as limits of flats $\{ E_i \}_{i \in \mathcal N}$ where $E_{i} \in \mathcal F_{i}$.

\begin{proposition} [See \cite{KL}, Proposition 4.3] \label{XinfProps}
The space $X_\omega$ satisfies the following two properties:
\begin{enumerate}
\item[\rm(F1)] Every non-degenerate open triangle in $X_\omega$ is contained in a flat $E \in \mathcal F_\omega$; and
\item[\rm(F2)] Any two flats in $\mathcal F_\omega$ intersect in at most a point.
\end{enumerate}
\end{proposition}
\begin{proof}
Let $\Delta = \Delta(x,y,z)$ be an open triangle in $X_\omega$.  Then $\Delta$ can be obtained as a limit of triangles $\Delta_{i}$, where $\Delta_{i} = \Delta (x_{i},y_{i},z_{i})$ is a triangle in $X_{i}$.  The triangle $\Delta_{i}$  may be identified with a triangle $\Delta'_{i}$ in $X$ (since $X$ and $X_n$ are the same set with different metrics).

For $\omega$-almost all $i$, the triangle $\Delta'_{i}$ is not $\delta$-thin, for otherwise the limit would not be a non-degenerate open triangle.  Therefore, $\Delta'_{i}$ is $\delta$-thin relative to some flat $E_{i} \in \mathcal F_X$.  Consider a point $w \in [x,y] \smallsetminus \{ x,y \}$.   The point $w \in X_\infty$ corresponds to a sequence of points $\{ w_{i} \}$.  Now $d_{\omega}(w,[y,z]) > 0$ and $d_{\omega}(w,[x,z])$, so for $\omega$-almost all $i$ the point $w_{i}$ is not contained in the $\delta$-neighbourhood of $[y_{i},z_{i}]$ or the $\delta$-neighbourhood of $[x_{i},z_{i}]$.  Therefore, for $\omega$-almost all $i$ the point $w_{i}$ is contained in the $\delta$-neighbourhood of $E_{i}$.  Let $u_{i}$ be a point in $E_{i}$ within $\delta$ of $w_{i}$.  It is clear that the sequences $\{ w_{i} \}$ and $\{ u_{i} \}$ have the same limit, namely $w$ (although $u_i$ is only defined for $\omega$-almost all $i$).  Therefore, $w$ is contained in the limit of the flats $\{ E_{i} \}$.  This proves Property (F1).

Now suppose that the flats $E, E' \in \mathcal F_\infty$ intersect in more than one point.  Let $x, y \in \hat{E}_1 \cap \hat{E}_2$ be distinct.  By Property (F1), there is a sequence of flats $\{ E_{i} \}$ which approximate $E$.  Let $u, w \in [x,y] \smallsetminus \{ x,y \}$ be arbitrary ($u \neq w$) and let $\{ x_i \}, \{ y_i \}, \{ u_i \}, \{ w_i \} \subseteq E_{i}$ be sequences of points representing $x$, $y$, $u$ and $w$, respectively.

Let $z \in E'$ be arbitrary so that $\Delta (x,y,z)$ is a non-degenerate triangle, and let $\{ z_i \in X_{i} \}$ be a sequence of points representing $z$.  Since the triangle $\Delta (x,y,z)$ is an open triangle, there is a sequence of flats $\{ E'_{i} \}$ whose limit contains $\Delta (x,y,z)$.

For $\omega$-almost all $i$, neither $u_i$ nor $w_i$ is contained in the $\delta$-neighbourhood of $[x_i,z_i] \cup [y_i,z_i]$, and so are contained in the $\delta$-neighbourhood of $E'_{i}$.  Therefore the points $u_i$ and $w_i$ are each contained in the $\delta$-neighbourhoods of both $E_{i}$ and $E'_{i}$.  However $u$ and $w$ are distinct, so for $\omega$-almost all $i$ the points $u_i$ and $w_i$ are at least $\phi(\delta)$ apart, which implies that $E_{i} = E'_{i}$ for $\omega$-almost all $i$.  Therefore, the triangle $\Delta(x,y,z)$ is contained in $E$.  Since $z$ was arbitrary, $E' \subseteq E$, and a symmetric argument shows that the two flats are equal.
\end{proof}

Using only the properties (F1) and (F2) from the conclusion of Proposition \ref{XinfProps} above, Kapovich and Leeb proved the following two results.

\begin{lemma} {\rm\cite[Lemma 4.4]{KL}}\label{ConstantProject}\qua
Let $E \in \mathcal F_\omega$ be a flat in $X_\omega$ and  let $\pi_E \co X_\omega \to E$ be the closest-point projection map.  Let $\gamma \co [0,1] \to X_\omega \smallsetminus E$ be a curve in the complement of $E$.  Then $\pi_E \circ \gamma \co [0,1] \to E$ is constant.
\end{lemma}

\begin{lemma} {\rm\cite[Lemma 4.5]{KL}}\label{NoLoops}\qua
Every embedded loop in $Y$ is contained in a flat $E \in \mathcal F_\omega$.
\end{lemma}

\subsection{The equivariant Gromov topology} \label{GromovSubsection}

From the sequence of homomorphisms $h_n \co G \to \Gamma$, we have constructed a space $X_{\omega}$, a basepoint $x_\omega$ and an isometric action of $G$ on $X_\omega$ with no global fixed point.    Let $X_\infty$ be the convex hull of the set $G . x_\omega$, and let $\mathcal C_\infty$ be the union of the geodesics $[x_\omega, g . x_\omega]$, along with the flats $E \in \Finf$ which contain some non-degenerate open triangle contained in a triangle $\Delta (g_1 . x_\omega, g_2 . x_\omega, g_3 . x_\omega)$, for $g_1,g_2,g_3 \in G$.  Certainly $X_\infty \subseteq \mathcal C_\infty$.  The set $\mathcal C_\infty$, and hence also $X_\infty$, is separable.

Note that $\mathcal C_\infty$ is a CAT$(0)$ space and that Proposition \ref{XinfProps}  and Lemmas \ref{ConstantProject} and \ref{NoLoops} hold for $\mathcal C_\infty$ also.  The action of $G$ on $X_\omega$ leaves $\mathcal C_\infty$ invariant, so there is an isometric action of $G$ on $\mathcal C_\infty$.  Since $C_\infty \subseteq X_\omega$, Lemma \ref{NoXinfFixedPt} implies:

\begin{lemma} \label{NoFixCinf}
There is no global fixed point for the action of $G$ on $\mathcal C_\infty$.
\end{lemma}

We have chosen to consider the space $\mathcal C_\infty$ rather than $X_\infty$ so that if some flat in $X_\omega$ intersects our subspace in a set containing a non-degenerate open triangle then the entire flat containing this triangle is contained in the subspace.

Suppose that $\{ (Y_n,\lambda_n) \}_{n=1}^\infty$ and $(Y,\lambda)$ are pairs consisting of metric spaces, together with actions $\lambda_n \co G \to \text{Isom}(Y_n)$, $\lambda \co G \to \text{Isom}(Y)$.  Recall (cf. \cite[\S 3.4, p. 16]{BF3}) that $(Y_n, \lambda_n ) \to (Y,\lambda)$ in the {\em $G$-equivariant Gromov topology} if and only if:  for any finite subset $K$ of $Y$, any $\epsilon > 0$ and any finite subset $P$ of $G$, for sufficiently large $n$, there are subsets $K_n$ of $Y_n$ and bijections $\rho_n \co K_n \to K$ such that for all $s_n, t_n \in K_n$ and all $g_1, g_2 \in P$ we have
\[	\left| d_{Y}(\lambda(g_1) . \rho_n (s_n) , \lambda(g_2) . \rho_n (t_n)) - d_{Y_n} ( \lambda_n(g_1) . s_n ,  \lambda_n(g_2) . t_n ) \right| < \epsilon .	\]

To a homomorphism $h \co G \to \Gamma$, we naturally associate a pair $(X_h,\lambda_h)$ as follows: let $X_h$ be the convex hull in $X$ of $G.x$ (where $x$ is the basepoint of $X$), endowed with the metric $\frac{1}{\mu_h}d_X$; and let $\lambda_h = \iota \circ h$, where $\iota\co \Gamma \to \text{Isom}(X)$ is the fixed homomorphism.

\begin{lemma} \label{GromovTop}
Let $\Gamma$, $X$, $G$ and $\{ h_n \co G \to \Gamma \}$ be as described above.  Let $X_\omega$ be the asymptotic cone of $X$, and $\mathcal C_\infty$ be as described above.  Let $\lambda_\infty \co G \to \text{Isom}(\mathcal C_\infty)$ denote the action of $G$ on $\mathcal C_\infty$ and $(\mathcal C_\infty, \lambda_\infty)$ the associated pair.

There exists a subsequence $\{ f_i \} \subseteq \{ h_i \}$ so that the elements $(X_{f_i},\lambda_{f_i})$ converge to $(\mathcal C_\infty, \lambda_\infty)$ in the $G$-equivariant Gromov topology.
\end{lemma} 
\begin{proof}
Since $\mathcal C_\infty$ is separable, there is a countable dense subset of $\mathcal C_\infty$, $S$ say.  Let $S_1 \subset S_2 \subset \ldots$ be a collection of finite sets whose union is $S$.  

Let $\{ 1 \} =Q_1 \subset Q_2 \subset \ldots \subset G$ be an exhaustion of $G$ by finite subsets.  Define $J_i$ to be the collection of $i \in \mathbb N$ so that $(X_{h_i}, \lambda_{i})$ is a $\frac{1}{i}$-approximation for $Q_i$ and $S_i$.  By the definition of asymptotic cone, $\omega (J_i) = 1$, and in particular each $J_i$ is infinite.

Let $n_1$ be the least element of $J_1$, and let $f_1 = h_{n_1}$.  Inductively, define $n_k$ to be the least element of $J_k$ which is not contained in $\{ n_1, \ldots , n_{k-1} \}$, and define $f_k = h_{n_k}$.

It is straightforward to see that the sequence $\{ f_i \}$ satisfies the conclusion of the lemma.
\end{proof}

The above result can be interpreted as a compactification of a certain space of metric spaces equipped with $G$-actions.  This is the `compactification' referred to in the title of this paper.

\begin{convention} \label{ConvergentSubseq}
For the remainder of the paper, we will assume that we started with the sequence $\{ f_i \co G \to \Gamma \}$ found in Lemma \ref{GromovTop} above.  This will allow us to speak of `all but finitely many $n$' instead of `$\omega$-almost all $n$'.
\end{convention}

To make the use of Convention \ref{ConvergentSubseq} more transparent, when using this convention we speak of the homomorphisms $f_i$, rather than $h_i$.  However, we still write $\lambda_i$ for the action of $G$ on $X$ induced by $f_i$, and we write $\mu_i$ for $\mu_{f_i}$ and $X_i$ for $X$ endowed with the metric $d_{X_i} := \frac{1}{\mu_i} d_X = \frac{1}{\mu_{f_i}} d_X$.

Let $\mathcal F_\infty$ be the set of flats in $\mathcal C_\infty$.  

\begin{corollary} \label{Flats}
Under Convention \ref{ConvergentSubseq}, for each $E \in \mathcal F_\infty$ there is a sequence $\{ E_i \subset X_i \}$ so that $E_i \to E$ in the $G$-equivariant Gromov topology.
\end{corollary}
\begin{proof}
This follows from the proofs of Proposition \ref{XinfProps} and Lemma \ref{GromovTop}.
\end{proof}

\subsection{The action of $G$ on $X_{\omega}$} \label{Action}

\begin{lemma} \label{FlatInv}
Let $E \in \mathcal F_\infty$ be a flat which is a limit of the flats $\{ E_i \}$.  If $g \in G$ and $g . E = E$ then for all but finitely many $j$ we have $f_j(g) . E_j = E_j$.
\end{lemma}
\begin{proof}
Choose a non-degenerate triangle $\Delta(a,b,c)$ in $E$.  Let $\{ E_{i} \}$ be a sequence of flats from $X_{i}$ approximating $E$.  Let $\{ a_i \}, \{ b_i \}$ and $\{ c_i \}$ be sequences of points representing $a,b$ and $c$, respectively.

By the definition of the action of $G$ on $X_\omega$, the point $g.x$ is represented by the sequence $\{ {\lambda}_{i}(g,a_i) \}$.  The triangle $\Delta(g.a,g.b,g.c)$ is also a non-degenerate triangle in $E$ and at least one of the triangles $\Delta(a,b,g.a)$, $\Delta (a,c,g.a)$, $\Delta (b,c,g.a)$ is non-degenerate in $E$.  The argument from the proof of Proposition \ref{XinfProps} (along with Corollary \ref{Flats}) applied to this non-degenerate triangle shows that for all but finitely many $i$ the point ${\lambda}_{i}(g,a_i)$ is $\delta$-close to the flat $E_{i}$.  Similarly,  for all but finitely many $i$ the point ${\lambda}_i(g,b_i)$ is $\delta$-close to $E_{i}$.  Since $E_i \in \mathcal F_X$, so is the flat $f_{i}(g) . E_{i}$, by Proposition \ref{FInv}.  However, 
\[	d_X({\lambda}_{i}(g,a_i), {\lambda}_i (g,b_i)) = d_X(a_i,b_i) ,		\]
which is greater than $\phi(\delta)$ for all but finitely many $i$. Therefore, by the definition of the function $\phi$, for all but finitely many $i$ the flats $E_i$ and $f_i(g) . E_{i}$ are the same, as required.
\end{proof}

\begin{lemma} \label{NoElliptic}
Suppose that $\Gamma$ is a group acting properly and cocompactly on a CAT$(0)$ space $X$ with isolated flats, and suppose that this action is toral.  Let $G$ and $X_\omega$ be as above.
Suppose that $g \in G$ leaves a flat $E$ in $\mathcal C_\infty$ invariant as a set.  Then $g$ acts by translation on $E$.
\end{lemma}
\begin{proof}
Suppose that $g$ acts nontrivially on $E$, but not as a translation.  Then there are $y, z \in E$ which are moved different distances by $g$ (suppose that $y$ is moved further than $z$ by $g$).  Let $\{ E_{i} \}$ be a sequence of flats in $X_{i}$ which converge to $E$.  Since $g$ maps $E$ to itself, for all but finitely many $i$, we have $f_i(g) . E_i = E_i$, by Lemma \ref{FlatInv}.

Let $\{ z_i \}, \{ y_i \} \subseteq E_{i}$ be sequences of points in representing $z$ and $y$, respectively.  Suppose that $d_{X_\omega}(y,g.y) - d_{X_\omega}(z,g.z) = \epsilon > 0$.  Choose large $i$ so that the points $z_i$, $y_i$, ${\lambda}_{i}(g , z_i)$ and ${\lambda}_{i}(g . y_i)$ satisfy
\begin{eqnarray*}
|d_{X_n}({\lambda}_{i}(g , z_i), z_i) - d_X(g.z,z)| & < & \frac{\epsilon}{3}; \mbox{ and }\\
|d_{X_n}({\lambda}_{i}(g , y_i), y_i) - d_X(g.y,y)| & < & \frac{\epsilon}{3}.
\end{eqnarray*}
Since $\Gamma$ is toral, the action of $g$ on $X_n$ via ${\lambda}_{i}$ is by (possibly trivial) translations.  Therefore,
\[	d_{X_n}({\lambda}_{i}(g, z_i), z_i) = d_{X_n}({\lambda}_{i}(g, y_i), y_i) .	\]
However, $d_X(g.y,y) - d_X(g.z,z) = \epsilon > 0$ and we have a contradiction.
\end{proof}

\begin{remark}
As we shall see in Example \ref{NontoralEx} below, Lemma \ref{NoElliptic} does not hold when $\Gamma$ is a non-toral CAT$(0)$ group with isolated flats.
\end{remark}

\subsection{Algebraic $\Gamma$-limit groups}

\begin{definition} (cf.\ \cite{Sela1}, Definition 1.2)\label{Kinf} \qua
Define the normal subgroup $K_\infty$ of $G$ to be the kernel of the action of $G$ on $\mathcal C_\infty$:
\[ K_\infty  = \{ g \in G\ |\ \forall y \in \mathcal C_\infty , \ g(y) = y \} . \]
The {\em strict $\Gamma$-limit group} is $L_\infty = G/K_\infty$. Let $\eta \co G \to L_\infty$ be the natural quotient map.

A {\em $\Gamma$-limit group} is a group which is either a strict $\Gamma$-limit group or a finitely generated subgroup of $\Gamma$.
\end{definition}

Recall the following (see \cite[Definition 1.5]{BF3}).

\begin{definition}
Let $G$ and $\Xi$ be finitely generated groups.  A sequence $\{ f_i\} \subseteq \text{Hom}(G, \Xi)$ is {\em stable} if, for all $g \in G$, the sequence $\{ f_i(g) \}$ is eventually always $1$ or eventually never $1$.  

For any sequence $\{ f_i \co G \to \Xi \}$ of homomorphisms, the {\em stable kernel} of $\{ f_i \}$, denoted $\SK (f_i)$, is
\[	\{ g\in G\ |\ f_i(g) = 1\ \mbox{ for all but finitely many $i$} \}.	\]
\end{definition}

\begin{definition}
An {\em algebraic $\Gamma$-limit group} is the quotient $G/\SK (h_i)$, where $\{ h_i \co G \to \Gamma \}$ is a stable sequence of homomorphisms.
\end{definition}

In the case that $\Gamma$ is a free group (acting on its Cayley graph),  Bestvina and Feighn  \cite{BF3} define {\em limit groups} to be those groups of the form $G/\SK f_i$, where $\{ f_i \}$ is a stable sequence in $\text{Hom}(G,\Gamma)$.  When $\Gamma$ is a free group, this leads to the same class of groups as the geometric definition analogous to Definition \ref{Kinf} above  (see \cite{Sela1}; this is also true when $\Gamma$ is a torsion-free hyperbolic group, see \cite{SelaHyp}).  When $\Gamma$ is a torsion-free CAT$(0)$ group with isolated flats we may have torsion in $G/K_\infty$, but $G/\SK$ is always torsion-free.  However, for any stable sequence $\{ f_i \}$, we always have $\SK f_i \subseteq K_\infty$.  Torsion in $G/ K_\infty$ can only occur when $\mathcal C_\infty$ is a single flat, in which case $f_i(G)$ is virtually abelian for almost all $i$.

In Section \ref{TreeSection} below, when $\Gamma$ is a torsion-free toral CAT$(0)$ group with isolated flats we use the action of $G$ on $\mathcal C_\infty$ to construct an action of $G$ on an $\R$-tree $T$.  In this case, the class of $\Gamma$-limit groups and algebraic $\Gamma$-limit groups coincides.  It will be this fact that allows us in \cite{CWIFShort} to prove many results about torsion-free toral CAT$(0)$ groups with isolated flats in analogy to Sela's results about free groups and torsion-free hyperbolic groups.

The following are elementary.

\begin{lemma} \label{SKinKinf}
Suppose that $\{ f_n \co G \to \Gamma \}$ gives rise to the action of $G$ on $\mathcal C_\infty$ as in the previous section.  Then $\SK ( f_n ) \subseteq K_\infty$.
\end{lemma}

\begin{lemma} \label{Limitfg}
Let $L$ be a $\Gamma$-limit group.  Then $L$ is finitely generated.
\end{lemma}

\subsection{A non-toral example} \label{NonToral}

In this paragraph we consider an example of the above construction in the case that $\Gamma$ is a torsion-free {\em non}-toral CAT$(0)$ group with isolated flats.

For any torsion-free group acting an properly and cocompactly on a CAT$(0)$ space $X$ with isolated flats, and for any maximal flat $E \in \Finf$, we know that $H:= \text{Stab}(E)$ is a torsion-free, proper cocompact lattice in $\mathbb R^n$, for some $n$.  Hence, by Bieberbach's theorem, $H$ has a free abelian group of finite index.  We have an exact sequence
\[	1 \to \Z^n \to H \to A \to 1,	\]
where $A$ is a finite subgroup of $O(n)$.  Each element $g \in H$ acts on $\R^n$ as
\[	g(v) = r_g(v) + t_g,	\]
where $r_g \in A \subset O(n)$ and $t_v \in \R^n$.  The homomorphism $H \to A$ is given by $g \to r_g$.

\begin{example} \label{NontoralEx}
Let $G_1$ be a non-abelian torsion-free crystallographic group as above, with the exact sequence
\[	1 \to \Z^n \to G_1 \to A \to 1,	\]
where $A$ is a nontrivial finite group and let $\Gamma = G_1 \times \mathbb Z$.  Let $w$ be the generator of the $\mathbb Z$ factor of  $\Gamma$.  Clearly the group $\Gamma$ acts properly and cocompactly by isometries on a CAT$(0)$ space with isolated flats.

Let $\Gamma$ be generated by $\{ g_1, \ldots , g_k , w \}$, where $\{ g_1, \ldots , g_k \}$ is a generating set for $G_1$ and let $F_{k+1}$ be the free group of rank $k+1$ with basis $\{ x_1 , \ldots , x_{k+1} \}$.  For $n \ge 1$, define the homomorphism $\phi_n \co F_{k+1} \to G_3$ by $\phi_n(x_i) = g_i$ for $1 \le i \le k$, and $\phi_n( x_{k+1}) = w^n$.  All of the kernels of $\phi_n$ are identical, so the algebraic $\Gamma$-limit group is $F / \SK (\phi_n) \cong \Gamma$.

In this case $\mathcal C_\infty = X_\infty = X_\omega = \mathbb R^{n+1}$.
In the geometric $\Gamma$-limit group, the $\mathbb Z^n$ in $G_1$ acts trivially, but the elements not in $\mathbb Z^n$ act like the corresponding element of $A$.  The element $w$ acts nontrivially by translation, and the $\Gamma$-limit group is isomorphic to $A \times \Z$, which is not torsion-free.
\end{example}

\section{The $\R$-tree $T$} \label{TreeSection}

For the remainder of the paper, we suppose that $\Gamma$ is a torsion-free {\em toral} CAT$(0)$ group with isolated flats.  In this section we extract an $\R$-tree $T$ from the space $X_\omega$ and an isometric action of $G$ on $T$.  The idea is to remove the flats in $X_\omega$ in order to obtain an $\R$-tree.  We replace the flats with lines.

\subsection{Constructing the $\R$-tree} \label{ConstructT}

Suppose that $\Gamma$ is a torsion-free toral CAT$(0)$ group with isolated flats acting on the space $X$, and that the sequence of homomorphisms $h_n \co G \to \Gamma$ gives rise to the limiting space $X_\omega$, as in the previous section.  Let $\mathcal C_\infty$ be the collection of geodesics and flats as described in Subsection \ref{Action}, and let $\mathcal F_\infty$ be the set of flats in $\mathcal C_\infty$.  

Suppose that $E \in \mathcal F_\infty$.  By Proposition \ref{XinfProps}, for any $g \in G$, exactly one of the following holds: (i) $g.E = E$; (ii) $| g.E \cap E| = 1$; or (iii) $g.E \cap E = \emptyset$.  By Lemmas \ref{FlatInv} and \ref{NoElliptic}, the action of $\text{Stab}(E)$ on $E$ is as a finitely generated free abelian group, acting by translations on $E$.

Let $\mathcal D_E$ be the set of directions of the translations of $E$ by elements of $\text{Stab}(E)$.  

For each element $g \in G \smallsetminus \text{Stab}(E)$, let $l_g(E)$ be the (unique) point where any geodesic from a point in $E$ to a point in $g.E$ leaves $E$, and let $\mathcal L_E$ be the set of all $l_g(E) \subset E$.  Note that if $g.E \cap E$ is nonempty (and $g \not\in \text{Stab}(E)$) then $g.E \cap E = \{ l_g(E) \}$.

Since $G$ is finitely generated, and hence countable, both sets $\mathcal D_E$ and $\mathcal L_E$ are countable.  
Given a (straight) line $p \subset E$, let $\chi_E^p$ be the projection from $E$ to $p$.  Since $\mathcal L_E$ is countable, there are only countably many points in $\chi_E^p(\mathcal L_E)$.  Therefore, there is a line $p_E \subseteq E$ such that
\begin{enumerate}
\item the direction of $p_E$ is not orthogonal to a direction in $\mathcal D_E$;
\item  if $x$ and $y$ are distinct points in $\mathcal L_E$, then $\chi_E^{p_E}(x) \neq \chi_E^{p_E}(y)$;
\end{enumerate}

Project $E$ onto $p_E$ using $\chi_E^{p_E}$.  The action of $\text{Stab}(E)$ on $p_E$ is defined in the obvious way (using projection) -- this is an action since the action of $\text{Stab}(E)$ on $E$ is by translations.  Connect $\mathcal C_\infty \smallsetminus E$ to $p_E$ in the obvious way -- this uses the following observation which follows immediately from Lemma \ref{ConstantProject}.

\begin{observation}
Suppose $S$ is a component of $\mathcal C_\infty \smallsetminus E$.  Then there is a point $x_S \in E$ so that $S$ is a component of $\mathcal C_\infty \smallsetminus \{ x_E \}$.
\end{observation}

Glue such a component $S$ to $p_E$ at the point $\chi_E^{p_E} ( x_S)$.

Perform this projecting and gluing construction in an equivariant way for all flats $E \subseteq \mathcal C_\infty$ -- so that for all $E \subseteq \mathcal C_\infty$ and all $g \in G$ the lines $p_{g . E}$ and $g . p_E$ have the same direction (this is possible since the action of $\text{Stab}(E)$ on $E$ is by translations, so doesn't change directions).

Having done this for all flats $E \subseteq \mathcal C_\infty$, we arrive at a space $T$, which is endowed with the (obvious) path metric.

The action of $G$ on $T$ is defined in the obvious way from the action of $G$ on $X_\omega$.  This action is clearly by isometries.

The space $T$ has a distinguished set of geodesic lines, namely those of the form $\chi_E^{p_E}(E)$, for $E \in F_\infty$.  Denote the set of such geodesic lines by $\mathbb P$.

\begin{lemma} \label{NoFixedPtOnT}
$T$ is an $\R$-tree and there is an action of $G$ on $T$ by isometries without global fixed points.
\end{lemma}
\begin{proof}
That $T$ is an $\R$-tree is obvious, since there are no embedded loops.  We have already noted that there is an isometric action of $G$ on $T$.

Finally, suppose that there is a fixed point $y$ for the action of $G$ on $T$.  If $y$ is not contained in some geodesic line in $\mathbb P$, then $y$ would correspond to a fixed point for the action of $G$ on $X_\omega$, and there are no such fixed points, by Lemma \ref{NoFixCinf}.

Thus $y$ is contained in some geodesic line $p_E \in \mathbb P$, corresponding to the flat $E \in \F_\infty$.  Let $g \in G$.  If $g$ does not fix $p_E$ then it takes $p_E$ to some line $p_{E'}$, and $g$ takes $E$ to $E'$ in $X_\omega$, fixing the point of intersection.  Suppose that $g_1$ and $g_2$ are elements of $G$ which fix $y$ but not $p_E$.  Then let $\alpha \in X_\omega$ be the point of intersection of $E$ and $g_1 . E$ and let $\beta$ be the point of intersection of $E$ and $g_2 . E$.  Then $\alpha$ and $\beta$ are both in $\mathcal L_E$ and $\chi_E^{p_E}(\alpha) = \chi_E^{p_E}(\beta)$, so by the choice of $p_E$ we must have $\alpha = \beta$.  Therefore, there is a point $\alpha \in E$ so that all elements $g \in G$ which do not fix $p_E \subseteq T$ fix $\alpha \in X_\omega$.

If $g$ does leave $p_E$ invariant then it fixes $E$ as a set, and so acts by translations on $E$, and hence by translations on $p_E$.  Therefore $g$ fixes $p_E$ pointwise, and the direction of translation of $g$ on $E$ is orthogonal to the direction of $p_E$.  By the choice of the direction of $p_E$ above, this means that $g$ acts trivially on $E$, and in particular fixes the point $\alpha$ found above.  Thus $\alpha$ is a global fixed point for the action of $G$ on $\mathcal C_\infty$, contradicting Lemma \ref{NoFixCinf}.
\end{proof}

\begin{remark}
Since $K_\infty \le G$ acts trivially on $\mathcal C_\infty$, it also acts trivially on $T$ and the action of $G$ on $T$ induces an isometric action of $L_\infty$ on $T$.
\end{remark}

\subsection{The actions of $G$ and $L_\infty$ on $T$}

The following theorem is the main technical result of this paper, and the remainder of this section is devoted to its proof.

Let $G$ be a finitely generated group, $\Gamma$ a torsion-free toral CAT$(0)$ group with isolated flats and $\{ h_i \co G \to \Gamma \}$ a sequence of homomorphisms, no two of which differ only by conjugation in $\Gamma$.  Let $X_\omega$, $\mathcal C_\infty$ and $T$ be as in Section \ref{Limits} and Subsection \ref{ConstructT} above.  Let $\{ f_i \co G \to \Gamma \}$ be the subsequence of $\{ h _i \}$ as in the conclusion of Lemma \ref{GromovTop}.  Let $K_\infty$ be the kernel of the action of $G$ on $\mathcal C_\infty$ and let $L_\infty = G / K_\infty$ be the associated strict $\Gamma$-limit group.

\begin{theorem} [Compare \cite{Sela1}, Lemma 1.3] \label{LinfProps}  In the above situation, the following properties hold.
\begin{enumerate}
\item Suppose that $[A,B]$ is a non-degenerate segment in $T$.  Then $\text{Fix}_{L_\infty}[A,B]$ is an abelian subgroup of $L_\infty$; \label{SegStabAb}
\item If $T$ is isometric to a real line then for all sufficiently large $n$ the group $f_n(G)$ is free abelian.  Furthermore in this case $L_\infty$ is free abelian. \label{Tline}
\item If $g \in G$ fixes a tripod in $T$ pointwise then $g \in \SK (f_i)$; \label{tripod}
\item Let $[y_1,y_2] \subset [y_3,y_4]$ be a pair of non-degenerate segments of $T$ and assume that the stabiliser $\text{Fix}([y_3,y_4])$ of $[y_3,y_4]$ in $L_\infty$ is non-trivial.  Then
\[	\text{Fix}([y_1,y_2]) = \text{Fix}([y_3,y_4]) .		\]
In particular, the action of $L_\infty$ on the $\R$-tree $T$ is stable. \label{Stable}
\item Let $g \in G$ be an element which does not belong to $K_\infty$.  Then for all but finitely many $n$ we have $g \not\in \text{ker}(f_n)$; \label{SeqStable}
\item $L_\infty$ is torsion-free; \label{tf}
\item If  $T$ is not isometric to a real line then $\{ f_i \}$ is a stable sequence of homomorphisms. \label{fiStable}
\end{enumerate}
\end{theorem}

We prove Theorem \ref{LinfProps} in a number of steps.

First, we prove \ref{LinfProps}(\ref{SegStabAb}).  Suppose that $[A,B] \subseteq T$ is a non-degenerate segment with a nontrivial stabiliser.  If there is a line $p_E \in \mathbb P$ such that $[A,B] \cap p_E$ contains more than one point, then any elements  $g_1,g_2 \in \text{Fix}([A,B])$ fix $p_E$ and hence fix $E \in \mathcal F_\infty$.  Therefore, by Lemma \ref{FlatInv} for all but finitely many $i$ the elements $g_1$ and $g_2$ fix the flat $E_i \in X_i$, where $\{ E_i \} \to E$.  The stabiliser of $E_i$ is free abelian, so $[ g_1, g_2] \in \text{ker}(f_i)$.   Thus $[g_1,g_2] \in \SK (f_i)$.  By Lemma \ref{SKinKinf}, $\SK (f_i) \subseteq K_\infty$, so $[g_1, g_2] \in K_\infty$.  Hence, in this case the stabiliser in $L_\infty$ of $[A,B]$ is abelian.

Suppose therefore that there is no $p_E \in \mathbb P$ which intersects $[A,B]$ in more than a single point.  In particular, $A$ and $B$ are not both contained in $p_E$ for any $p_E \in \mathbb P$. 

In fact, we need something stronger than this.  First we prove the following:

\begin{lemma} \label{NearFlat}
Suppose that $\alpha, \beta \in X_n$ and $g \in G$ are such that there is a segment of length at least $\SegCons + 4\max \{ d_{X}(g.\alpha, \alpha), d_{X}(g. \beta, \beta) \}$ in $[\alpha, \beta]$ which is within $\delta$ of a flat $E \in \F_X$.  Then $g \in \text{Fix}(E)$.
\end{lemma}
\begin{proof}
In this lemma, all distances are measured with the metric $d_X$.  Let $L = \max \{ d_{X}(g.\alpha, \alpha), d_{X}(g. \beta, \beta) \}$.

Let $[\alpha_1, \beta_1]$ be the segment in $[\alpha,\beta]$ of length at least $\SegCons + 4L$ which is in the $\delta$-neighbourhood of $E$.  

Consider first the triangle $\Delta_1 = \Delta(\alpha, \beta, g . \beta)$.

If $\Delta_1$ is $\delta$-thin, then there is a segment $[\alpha_2,\beta_2] \in [\alpha, g . \beta]$ which has length at least $\SegCons + 3L - 2\delta$ and is within $\delta$ of $[\alpha_1,\beta_1]$.  Hence $[\alpha_2,\beta_2]$ is in the $2\delta$-neighbourhood of $E$.  Also, when $[\alpha, \beta]$ and $[\alpha, g . \beta]$ are both parametrised by arc length, there is an interval of time of length at least $\SegCons + 3L - 2\delta$ when $[\alpha_1,\beta_1]$ and $[\alpha_2,\beta_2]$ both occur.

Suppose then that $\Delta_1$ is $\delta$-thin relative to a flat $E' \neq E$.  Then a subsegment $[\alpha_1', \beta'] \subseteq [\alpha_1,\beta_1]$ of length at least $\SegCons +3L - \phi(\delta)$ does not intersect the $\delta$-neighbourhood of $E'$, and so must be contained in the $\delta$-neighbourhood of $[\alpha, g. \beta] \cup [\beta, g . \beta]$.  However, $d_X(\alpha,\beta) \ge \SegCons + 4L$, so there is an interval $[\alpha_2, \beta_2]$ in $[\alpha, g . \beta]$ of length at least $\SegCons + 3L - \phi(\delta)$ which is within $\delta$ of $[\alpha_1, \beta_1]$ and hence in the $2\delta$-neighbourhood of $E$.  Once again, when $[\alpha, \beta]$ and $[\alpha, g . \beta]$ are parametrised by arc length there is an interval of time of length at least $\SegCons + 3L - \phi(\delta)$ when both $[\alpha_1,\beta_1]$ and $[\alpha_2, \beta_2]$ occur.

Finally suppose that $\Delta_1$ is $\delta$-thin relative to $E$.  Then there is certainly a segment $[\alpha_2,\beta_2] \subseteq [\alpha, g. \beta]$ of length at least $\SegCons + 3L - 2\delta$ in the $\delta$-neighbourhood of $E$. Since 
\begin{eqnarray*}
d_X(\alpha, g . \beta) & \le & d_X(\alpha , \beta) + d_X(\beta, g . \beta)\\
& \le & d_X(\alpha, \beta) + L,
\end{eqnarray*}
then just as above there are subsegments $[\alpha_1',\beta_1']$ and $[\alpha_2' , \beta_2']$ (contained in $[\alpha_1,\beta_1]$ and $[\alpha_2,\beta_2]$ ,respectively) which occur at the same time for an interval of at least $\SegCons + 2L - 2\delta$.

In any of these cases, denote by $[\alpha_1, \beta_1]$ and $[\alpha_2, \beta_2]$ the intervals in $[\alpha , \beta]$ and $[\alpha , g . \beta]$ of length at least $\SegCons + 2L - 2\phi(\delta)$ which occur at the same time when $[\alpha, \beta]$ and $[\alpha, g . \beta]$ are parametrised by arc length and such that $[\alpha_1,\beta_1]$ and $[\alpha_2,\beta_2]$ are both contained in the $2\delta$-neighbourhood of $E$.  We now consider $\Delta_2 = \Delta(g . \alpha, \beta, g . \beta)$.  

Using the same arguments as those which found $[\alpha_2,\beta_2]$ as above, it is not difficult to find an interval $[\alpha_3,\beta_3] \subseteq [g.\alpha, g . \beta]$ which is of length at least $\SegCons - 4\phi(\delta)$ and occurs at the same time as $[\alpha_2,\beta_2]$ when $[\alpha, g . \beta]$ and $[g . \alpha, g . \beta]$ are parametrised by arc length.  The only wrinkle in this argument occurs when $\Delta_2$ is $\delta$-thin relative to $E$ and we may have to change $[\alpha_2,\beta_2]$ as in the third case above.  However, in this case we can find an appropriate $[\alpha_1,\beta_1]$.

Now, $[\alpha_3,\beta_3]$ is contained in the $4\delta$-neighbourhood of $E$.  Also, the time at which it occurs overlaps the time at which $[\alpha_1,\beta_1]$ occurs by at least $\SegCons - 4\phi(\delta)$.  Note that $\SegCons - 4\phi(\delta) > \phi(4\delta) + 1$.  Since $[\alpha_1,\beta_1] \subseteq [\alpha, \beta]$ is contained in the $\delta$-neighbourhood of $E$, and since $[\alpha_3, \beta_3] \subseteq [g . \alpha, g . \beta]$ occurs at the same time as $[\alpha_1, \beta_1]$, the interval $[\alpha_3, \beta_3]$ is contained in the $\delta$-neighbourhood of $g . E$.  Therefore, the $4\delta$-neighbourhoods of $E$ and $g . E$ intersect in a geodesic segment of length at least $\phi(4\delta) + 1$, which implies that $E = g . E$, so $g \in \text{Fix}(E)$, as required.
\end{proof}

Fix a finite subset $Q \in \text{Fix}_G([A,B])$.  Let the sequences $\{ A_k \}$ and $\{ B_k \}$ converge to $A$ and $B$, respectively.

Suppose that, for some $\epsilon > 0$, for all but finitely many $i$ there is a segment $[\alpha_{i},\beta_{i}] \subseteq [A_{i},B_{i}]$ for which (i) $d_{X_i}(\alpha_{i},\beta_{i}) \ge \epsilon$; and (ii) there is a flat $E_{i}$ so that $[\alpha_{i},\beta_{i}]$ is in the $\delta$-neighbourhood of $E_{i}$.  In this case, for all but finitely many $i$, the conditions of Lemma \ref{NearFlat} are satisfied for each $g \in Q$ and the segment $[\alpha_{i},\beta_{i}]$.  Therefore, $f_{i}(Q) \subseteq \text{Stab}_{\Gamma}(E_{i})$, so $\langle f_{i}(Q) \rangle$ is abelian for all but finitely many $i$.   In this case $\left\langle q K_\infty \ |\ q \in Q \right\rangle$ is certainly abelian.  Thus we can suppose that for all $\epsilon > 0$ there is no such segment $[ \alpha_i, \beta_i ]$.

Let $Q = \{ g_1 , \ldots , g_s \}$.  Let $[\rho, \sigma] \in [A, B]$ be the middle third, and let the sequences $\{ \rho_k \}$, $\{ \sigma_k \}$, $\{ {\lambda}_k(g_i, \rho_k) \}$ and $\{ {\lambda_k}(g_i, \sigma_k) \}$ converge to $\rho$, $\sigma$, $g_i . \rho$, and $g_i . \sigma$.  Consider the triangles $\Delta (\rho_k, \sigma_k, \lambda(g_i , \sigma_k))$ and $\Delta( \lambda(g_i , \rho_k), \sigma_k, \lambda(g_i , \sigma_k))$.  By the argument in Lemma \ref{NearFlat}, the argument in the above paragraph and the argument in \cite[Proposition 2.4]{Paulin}, for all but finitely many $k$ there is a segment $[\tau_k, \upsilon_k] \subseteq [\rho_k, \sigma_k]$ whose length (measured in $d_{X_k}$) is at least $\frac{1}{3}d_{X_k}(\rho_k,\sigma_k)$ and such that for each $g_i \in Q$, the image $g_i . [\tau_k,\upsilon_k]$ lies in the $4\delta$-neighbourhood of $[\rho_k,\sigma_k]$.  Note that $d_{X_k}(\tau_k,\upsilon_k) \ge \frac{1}{9}(A_k,B_k)$.

Now, since translations on a line commute, each element of the form $[g,g']$, $g,g' \in Q$ moves the midpoint of $[\tau_k,\upsilon_k]$ at most $16\delta$ (see \cite[Proposition 2.4]{Paulin}).

There is an absolute bound on the size of the ball of radius $16\delta$ around any point in $X$.  Let this bound be $D_1$.  Since the action of $\Gamma$ is free, we have $|\{ [f_{n}(g),f_{n}(g')] \ | g,g' \in Q \}| \le D_1$ irrespective of the size of $Q$.  

Let $x_1, x_2 \in \text{Fix}_G([A,B])$.  Either for all but finitely many $n$ there is a flat $E_n$ so that $x_1, x_2 \in \text{Fix}(E_n)$ or the above argument bounding the size of the set of commutators holds.  In the first case, $\langle f_{n}(x_1), f_{n}(x_2) \rangle$ is abelian for all but finitely many $n$.

Suppose then that the second case holds.  By the above argument, with $Q = \{ x_2, x_1, x_1^2, \ldots , x_1^{D_1+1} \}$, for all but finitely many $n$ we can find $1 \le s_1 < s_2 \le D_1 + 1$ for which $f_{n}([x_1^{s_1},x_2]) = f_{n}([x_1^{s_2},x_2])$.  Therefore $f_{n}(x_1^{s_1}x_2x_1^{-s_1}) = f_{n}(x_1^{s_2}x_2x_1^{-s_2})$, which implies that $f_{n}([x_1^{s_2-s_1},x_2]) = 1$.  Hence $f_{n}(x_1)^{s_2-s_1}$ commutes with $f_{n}(x_2)$.  Also, $\langle f_{n}(x_1)^{s_2-s_1} \rangle \subseteq \langle f_{n}(x_1)^{s_2-s_1} , f_{n}(x_2) \rangle^{f_{n}(x_1)}$, so by Proposition \ref{malnormal}, $f_{n}(x_1)$ is contained in the same maximal abelian subgroup as $\langle f_{n}(x_1)^{s_2 - s_1}, f_{n}(x_2) \rangle$, which is to say that $f_{n}(x_1)$ commutes with $f_{n}(x_2)$.

Therefore, in any case if $x_1, x_2 \in \text{Fix}_G([A,B])$ then for all but finitely many $n$ $f_{n}(x_1)$ commutes with $f_{n}(x_2)$.  Therefore $x_1 K_\infty$ commutes with $x_2 K_\infty$.  Since $x_1$ and $x_2$ were arbitrary, we have proved that the group $\text{Fix}_{L_\infty}([A,B])$ is abelian. This finishes the proof of \ref{LinfProps}(\ref{SegStabAb}).

\medskip

We now prove \ref{LinfProps}(\ref{Tline}).  Suppose that $T$ is isometric to a real line, so $L_\infty$ is a subgroup of $\text{Isom}(\R)$. 

Suppose first that $\mathcal C_\infty$ is not a single flat.  

Suppose that $k_1,k_2 \in G$ are arbitrary.  Let $H = \langle k_1,k_2 \rangle$ and $\overline{H}$ be the image of $\langle k_1,k_2 \rangle$ in $L_\infty$.  Then $\overline{H}$ is a $2$-generator subgroup of $\text{Isom}(\R)$ and so is one of the following (i) cyclic; (ii) infinite dihedral; or (iii) free abelian of rank $2$.  We first prove that $\overline{H}$ cannot be infinite dihedral, so that it is abelian.  

Suppose that $k \in G$ reverses the orientation of $T$ (which we are assuming is isometric to $\R$).  Let $a,b \in T$ be distinct points so that $k(a) = b$ and $k(b) = a$.  Approximate the segment $[a,b]$ by a segment $[a_i,b_i] \subset X_i$ for large $i$, and let $[c_i,d_i]$ be the middle third of $[a_i,b_i]$.  Since $\mathcal C_\infty$ is not a single flat, the segment $k.[c_i,d_i]$ lies within $2\delta$ of $[a_i,b_i]$, with orientation reversed (distances are being measured in $d_X$).  It is not hard to see that there must be a point $e_i \in [c_i,d_i]$ which is moved at most $2\delta$ by $k$.  Therefore $d_X(k^2 . e_i, e_i) \le 6\delta$.  In turn, this implies that $k^2$ moves each point on $[c_i,d_i]$ distance at most $10\delta$.  Repeating this argument (with a larger $i$) with the elements $k, k^3, \ldots , k^{2D_1+1}$ we find $1 \le i_1 < i_2 \le D_1$ so that $k^{2(2i_1+1)} . e_i = k^{2(2i_2+1)} . e_i$, which implies that $k^{2(2i_1+1)} k^{-2(2i_2+1)} \in \text{ker}(f_i)$ for large enough $i$.  This in turn implies, since $\Gamma$ is torsion-free, that $k \in \text{ker}(f_i)$ for all but finitely many $i$.  Therefore $k$ acts trivially on $T$, and so cannot reverse the orientation.  This proves that $\overline{H}$ is abelian, and being an orientation preserving subgroup of $\text{Isom}(\R)$, it is free abelian.

Suppose then that $k_1$ and $k_2$ act as translations on $T$, with translations lengths $\tau_1, \tau_2$, say.  Choose $\kappa > 100D_1 ( \max \{ \tau_1, \tau_2 \} + 1)$, and choose $a,b \in T$ distance $\kappa$ from each other.  For large enough $i$, the approximation $[a_i,b_i] \subset X_i$ of $[a,b]$ is such that for all $j \in \{ 1, \ldots , D_1 \}$, the elements $k_1^j$ and $k_2$ move the middle third $[c_i,d_i]$ of $[a_i,b_i]$ entirely within the $4\delta$ neighbourhood of $[a_i,b_i]$.  As in the proof of \ref{LinfProps}(\ref{SegStabAb}) above, the commutator $[k_1^j,k_2]$ moves the midpoint of $[c_i,d_i]$ a distance at most $16\delta$.  Therefore, there is $1 \le j_1 < j_2 \le D_1$ so that $f_i([k_1^{j_1},k_2]) = f_i([k_1^{j_2},k_2])$, which as above implies that $f_i([k_1,k_2]) = 1$ for large enough $i$.  

In this argument $k_1$ and $k_2$ were arbitrary, so letting $k_1$ and $k_2$ run over all pairs in a finite generating set $\{ g_1, \ldots ,g_k \}$ for $G$ we see that for all $i, j \in \{ 1 , \ldots , k \}$ the elements $f_{n}(g_i)$ and $f_{n}(g_j)$ commute for all but finitely many $n$.  Thus, $f_{n}(G)$ is abelian for all but finitely many $n$.  

We have also proved that $L_\infty$ is an orientation preserving subgroup of $\text{Isom}(\R)$, which is free abelian, as required.

Now suppose that $\mathcal C_\infty$ is a single flat.  By Lemma \ref{FlatInv}, for any $g \in G$ for all but finitely many $i$ the element $g$ fixes a flat $E_i \subseteq X_i$.  Take a finite generating set for $G$ and note that for all but finitely many $i$ each of the elements in this set fix the flat $E_i$.  Therefore, for all but finitely many $i$, $f_{i}(G) \subseteq \text{Fix}(E_i)$, which is free abelian.  This proves also that $L_\infty$ is abelian, and again the only abelian subgroups of $\text{Isom}(\R)$ which are not free abelian have a global fixed point.
This proves \ref{LinfProps}(\ref{Tline}).  

\medskip

We now prove \ref{LinfProps}(\ref{tripod}).  Let $T(A,B,C)$ be a tripod in $T$ and let $N$ be the valence three vertex in $T(A,B,C)$.  Suppose that $g \in G \smallsetminus \{ 1 \}$ stabilises $A, B$ and $C$ and therefore also $N$.  We prove that $g \in \text{ker}(f_n)$ for all but finitely many $n$. 

Let $K_0$ be the maximum number of elements of any orbit $\Gamma . y$ in any ball of radius $170\phi(16\phi(\phi(3\delta)))$ in $X$ (with metric $d_X$).  Such a $K_0$ exists because the action of $\Gamma$ on $X$ is proper and cocompact.

Suppose that some point $\alpha \in T(A,B,C)$ is contained in a line $p_E \in \mathbb P$.  Certainly not all of $T(A,B,C)$ is contained in $p_E$, so let $y \in T(A,B,C) \smallsetminus p_E$. Let $\alpha \in \mathcal C_\infty$ be a point corresponding to $y \in T$ and $\beta = \pi_E(\alpha)$.  Let $\{ E_i \}$ be a sequence of flats converging to $E$ (such a sequence exists by Corollary \ref{Flats}).  By Lemma \ref{FlatInv} $f_i(g) \in \text{Fix}(E_i)$ for all but finitely many $i$.

By fixing a sufficiently small $\epsilon$ and finding an $X_i$ which is an $\epsilon$ approximation for $Q = \{ 1, g , \ldots , g^{K_0+1} \}$ and $\{ \alpha, \beta \}$ we can ensure that for all $q \in Q$ the geodesic $[\alpha_i, {\lambda}_{i}(q, \alpha_i)]$ does not intersect the $4\delta$-neighbourhood of $E_i$, where $E$ is the limit of the flats $\{ E_i \}$.  Hence by Proposition \ref{Projection} $d_X(\pi_{E_i}(\alpha_i), \pi_{E_i}({\lambda}_{i}(q , \alpha_i))) \le 2\phi(3\delta)$, for all $q \in Q$.  However, $|Q| > K_0$ and there are no more than $K_0$ elements of $\Gamma . y$ in a ball of radius $2\phi(3\delta)$ in $X$, by the choice of $K_0$.  Note also that if $f_i(g)$ leaves $E_i$ invariant then $\pi_{E_i}({\lambda}_{i}(q, \alpha_i)) = {\lambda}_i(q, \pi_{E_i}(\alpha_i))$. Therefore, there exist $1 \le s_1 < s_2 \le K_0 + 1$ so that ${\lambda}_i (g^{s_1}, \pi_{E_i}(\alpha_i)) = {\lambda}_i(g^{s_2}, \pi_{E_i}(\alpha_i))$, which implies that $f_i(g^{s_2-s_1})$ fixes $\pi_{E_i}(\alpha_i)$.  Therefore, $f_i(g^{s_2-s_1}) = 1$, and since $\Gamma$ is torsion-free $g \in \text{ker}(f_{i})$ for all but finitely many $i$, as required.

Therefore, we assume for the moment that no point in the tripod $T(A,B,C)$ is contained in any line $p_E \in \mathbb P$.  In this case, $A,B$ and $C$ correspond to points $\A, \B$ and $\C$ in $\mathcal C_\infty$ for which $\Delta(\A,\B,\C)$ is a tripod in $\mathcal C_\infty$.   Let $\overline{N}$ be the valence three vertex in the tripod $T(\A,\B,\C)$.

\begin{figure}[ht!]\anchor{ApproxPic}
\begin{center}
\scalebox{1.5}{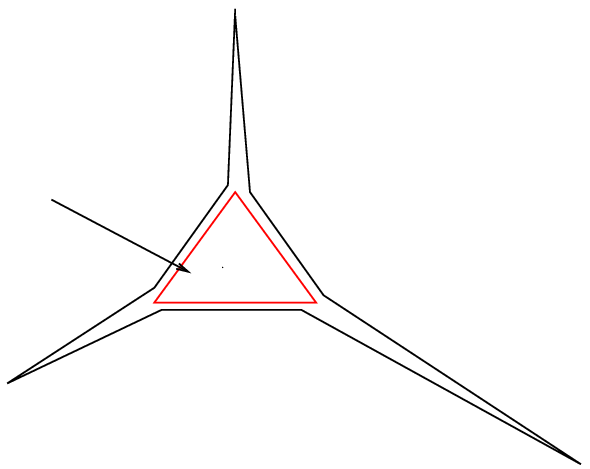}
\end{center}
\caption{Ensuring $\A, \B$ and $\C$ do not lie close to the flat $E$}
\label{ApproxPic}
\end{figure}

Let $\A', \B', \C'$ be the midpoints of $[\A,\overline{N}], [\B,\overline{N}]$ and $[\C,\overline{N}]$ respectively and let $S = \{ \A, \B, \C, \overline{N}, \A', \B', \C' \}$.  We also define the set $Q  = \{ 1, g , g^2, \ldots , g^{K_0 +1} \}$.  For varying $\epsilon$, we will consider those $X_i$ which are $\epsilon$ approximations for $Q$ and $S$.  Consider the triangle $\Delta = \Delta(\A_k,\B_k,\C_k)$ in $X_k$, an $\epsilon$-approximation for $Q$ and $S$.  Suppose that $\Delta$ is $\delta$-thin relative to a flat $E$.  If $\epsilon$ is small enough, then necessarily $\A_k$ is at least distance $\frac{1}{3}d_{X_i}(\A_k, \A_k')$ from $E$.  See \figref{ApproxPic}.

Note that since $\A_k$, $\B_k$ and $\C_k$ are not moved far by $q \in Q$, compared to the distances $d_{X_i}(\A_k,\A_k')$, etc., the same property is true for triangles such as $\Delta(\A_k, \B_k, q . \C_k)$.  

Fix $X_i$, an $\epsilon$-approximation for $Q$ and $S$ so that $\epsilon$ is `small enough' in the sense of the previous two paragraphs, and also $\epsilon < \frac{1}{100}$.

Consider the triangle $\Delta = \Delta(\A_i, \B_i, \C_i)$ in $X_i$.  Define the constant $\delta' = 17\phi(16\phi(\phi(3\delta)))$.

Suppose that $\Delta$ is not $\delta'$-thin, so it is $\delta$-thin relative to a unique flat $E \subset X_i$.  In this case $d_X(\pi_E(\A_i), \pi_E(\B_i)) \ge \delta' - 2\delta$.  Now, since $\epsilon < \frac{1}{100}$ and $d_X(\A_i, \pi_E(\A_i)) \ge \frac{1}{3}d_{X_i}(\A_k,\A_k')$, for large enough $i$ the geodesic $[\A_i, {\lambda}_i(q , \A_i)]$ avoids the $4\delta$-neighbourhood of $E$ for all $q \in Q$.  Therefore, for all $q \in Q$, $d_X(\pi_E(\A_i), \pi_E({\lambda}_i(q , \A_i))) \le 2 \phi(3\delta)$.  Fix $q \in Q$.  We now prove that $f_i(q)$ leaves $E$ invariant, and then as above we argue that $g \in \text{ker}(f_i)$.

To prove that $f_i(q)$ leaves $E$ invariant, we consider the triangle 
\[ \Delta' = \Delta({\lambda}_i(q , \A_i),  {\lambda}_i(q , \B_i) , {\lambda}_i(q , \C_i)),	\]
and prove that it is $16\phi(\phi(3\delta))$-thin relative to $E$.  Since it is also $\delta$-thin relative to $f_i(q) . E$ and since it is not $\delta'$-thin, we must have that $E = f_i(q) . E$, by an argument similar to that which proved Lemma \ref{UniqueFlat}.   

Let $\alpha_1 = {\lambda}_i(q,\A_i), \alpha_2 = \pi_E(\alpha_1), \beta_1 = {\lambda}_i(q, \B_i)$ and $\beta_2 = \pi_E(\beta_1)$.  Now, 
\begin{eqnarray*}
d_{X}(\alpha_2, \pi_E(\A_i)), d_{X}(\beta_2), \pi_E(\B_i)) & \le & 2\phi(3\delta), \mbox{ and }\\ d_{X}(\pi_E(\A_i), \pi_E(\B_i)) & \ge & \delta' - 2\delta,
\end{eqnarray*}
so we have
\[	d_{X}(\alpha_2, \beta_2) \ge 17\phi(16\phi(\phi(3\delta))) - 4\phi(3\delta) - 2\delta .	\]

Now, by Lemma \ref{Triangle} there exists $u,v$ in the $2\delta$-neighbourhood of $E$ so that $u \in [\alpha_1, \alpha_2]$ and $v \in [\alpha_1, \beta_2]$ and $d_{X}(u,v) \le \phi(\delta)$.  Now, since $\alpha_2 = \pi_E(\alpha_1)$, $d_{X}(u,\alpha_2) \le 2\delta$, and so $d_{X}(v,\alpha_2) \le \phi(\delta) + 2\delta)$.  Therefore, $[\alpha_2,\beta_2]$ and $[v,\beta_2]$ $3\phi(\delta)$-fellow travel and $\Delta(\alpha_1,\alpha_2,\beta_2)$ is $3\phi(\delta)$-thin.  Also $d_{X}(v,\beta_2) \ge \delta' - 7\phi(3\delta) - 2\delta$.

Now, consider the triangle $\Delta_1 = \Delta(\alpha_1,\beta_1,\beta_2)$.  Suppose it is $\delta$-thin relative to a flat $E' \neq E$.  Let $w_1$ be the point on $[\beta_2,\beta_1]$ which lies in the $4\phi(\delta)$-neighbourhood of $E$ furthest from $\beta_2$.  Then $w_1$ is not in the $\delta$-neighbourhood of $[v,\beta_2]$, and so is not in the $\delta$-neighbourhood of $[\alpha_1,\beta_2]$.  If $w_1$ is in the $\delta$-neighbourhood of $[\beta_1,\alpha_1]$ then $\Delta_1$ is $5\phi(\delta)$-thin.  Thus suppose that $w_1$ is in the $\delta$-neighbourhood of $E'$.  

Now, because (i) $\beta_2 = \pi_E(\beta_1)$; (ii) $d_X(v,\beta_2) \ge \delta' - 7\phi(3\delta) - 2\delta$; and (iii) $[v,\beta_2]$ is contained in the $3\phi(\delta)$-neighbourhood of $E$, either $\Delta_1$ is $(\phi(3\phi(\delta)) + 3\phi(\delta) + \delta)$-thin or there is a segment in $[v,\beta_2]$ of length at least $\phi(\phi(3\delta))$ which lie in the $\delta$-neighbourhood of $E'$.  However, this segment also lies in the $\phi(3\delta)$-neighbourhood of $E$, which is a contradiction.  Therefore, in any case either  $\Delta_1$ is $5\phi(\phi(3\delta))$-thin or $\Delta_1$ is $\delta$-thin relative to $E$.

Suppose that $\Delta_1$ is $\delta$-thin relative to $E$.  The geodesic $[\beta_2,\beta_1]$ intersects the $\delta$-neighbourhood of $E$ in a segment of length at most $\delta$ and so in this case $\Delta_1$ is $2\delta$-thin.

We have proved that $\Delta_1 = \Delta(\alpha_1,\beta_1,\beta_2)$ is $5\phi(\phi(3\delta))$-thin, and also that $\Delta(\alpha_1,\alpha_2,\beta_2)$ is $3\phi(\delta)$-thin.  Therefore, the geodesic $[\alpha_1,\beta_1]$ must $8\phi(\phi(3\delta))$-fellow travel the path $[\alpha_1,\alpha_2,\beta_2,\beta_1]$.

Similar arguments applied to the geodesic segments $[{\lambda}_i(q, \A_i), {\lambda}_i(q , \C_i)]$ and $[{\lambda}_i(q, \B_i), {\lambda}_i(q , \C_i)]$ show that the triangle $\Delta' = \Delta({\lambda}_i(q , \A_i),  {\lambda}_i(q , \B_i) , {\lambda}_i(q , \C_i))$ is $16\phi(\phi(3\delta))$-thin relative to $E$.  Since $\Delta'$ is not $\delta'$-thin, the argument from Lemma \ref{UniqueFlat} implies that $\Delta'$ is $16\phi(\phi(3\delta))$-thin relative to a unique flat.  Since $\Delta'$ is certainly $16(\phi(3\delta))$-thin relative to $f_i(q) . E$ we must have that $f_i(q) . E = E$, as required.  Now we know for all $q \in Q$ that $d_X({\lambda_i}(q, \pi_E(\A_i)), \pi_E(\A_i)) \le 2\phi(3\delta)$, which implies as above that $g \in \text{ker}(f_i)$.  

Therefore, we may assume that $\Delta$ is $\delta'$-thin.  Similar arguments to those above allow us to infer that for all $r_1, r_2, r_3 \in \{ 1 , q \}$, the triangle 
\[  \Delta(f_i(r_1) . \A_i, f_i(r_2) . \B_i, f_i(r_3) . \C_i) ,	\]
is $\delta'$-thin.  Now, by the Claim in the proof of \cite[Lemma 4.1]{RipsSelaGAFA}, the point $\overline{N}_i$ is moved by $f_i(q)$ at most $170\phi(16\phi(\phi(3\delta)))$.  Again in this case, we find $1 \le s_1 < s_2 \le K_0 +1$ so that $f_i(g^{s_2-s_1})$ fixes $\overline{N}_i$, which implies that $g \in \text{ker}(f_i)$.

We have proved that if $g \in G$ stabilises a tripod in $T$ then for all but finitely many $i$ $g \in \text{ker}(h_i)$.  This proves \ref{LinfProps}(\ref{tripod}).

\medskip

The proof of \ref{LinfProps}(\ref{Stable}) is identical to the of \cite[Proposition 4.2]{RipsSelaGAFA}, except that segment stabilisers are abelian, rather than cyclic.  However, all that is used in this proof is that segment stabilisers are abelian.

\medskip

We now prove \ref{LinfProps}(\ref{SeqStable}).  Suppose that $g \not\in K_\infty$.  Then certainly $g \not\in \text{ker}(f_{k})$ for all but finitely many $k$, by the choice of the sequence $\{ f_i \}$ in Lemma \ref{GromovTop}.

\medskip

We now prove \ref{LinfProps}(\ref{tf}).  If $T$ is isometric to a real line then by \ref{LinfProps}(\ref{Tline}) $L_\infty$ is finitely generated free abelian, and is certainly torsion-free.
  
Therefore suppose that $T$ is not isometric to a real line, that $g \in G$ and that $g^p \in K_\infty$.  Since $T$ is not isometric to a real line, $g^p$ stabilises a tripod, so $g^p \in \text{ker}(f_{k})$ for all but finitely many $k$.  However $\Gamma$ is torsion-free, so $g \in \text{ker}(f_{k})$ for all but finitely many $k$ and, by \ref{LinfProps}(\ref{SeqStable}), $g \in K_\infty$, as required.  This proves \ref{LinfProps}(\ref{tf}).

\medskip

Finally, we prove \ref{LinfProps}(\ref{fiStable}).  To see that $\{ f_i \}$ is a stable sequence of homomorphisms when $T$ is not isometric to a real line, suppose that $g \in G$.  If $g \not\in K_\infty$ then by \ref{LinfProps}.(\ref{SeqStable}) we have $g \not\in \text{ker}(f_i)$ for all but finitely many $i$.  If $g \in K_\infty$ then $g$ stabilises a tripod in $T$ and so by \ref{LinfProps}.(\ref{tripod}) $g \in \text{ker}(f_i)$ for all but finitely many $i$.  This proves that $\{ f_i \}$ is a stable sequence of homomorphisms.

This finally completes the proof of Theorem \ref{LinfProps}. \hfill $\square$

\section{$\Gamma$-limit groups and concluding musings} \label{conclusion}

\subsection{Various kinds of $\Gamma$-limit groups}

\begin{theorem}
Suppose that $\Gamma$ is a torsion-free toral CAT$(0)$ group with isolated flats.  The class of $\Gamma$-limit groups coincides with the class of algebraic $\Gamma$-limit groups.
\end{theorem}
\begin{proof} 
The class of abelian $\Gamma$-limit groups is exactly class of finitely generated free abelian groups.  It is easy to see that these are also the abelian algebraic $\Gamma$-limit groups.

Clearly a finitely generated subgroup of $\Gamma$ is an algebraic $\Gamma$-limit group.

Suppose then that $\{ h_i \co G \to \Gamma\}$  is a sequence of homomorphisms and $\{ f_i \}$ is the subsequence obtained from Lemma \ref{GromovTop}.  If the limiting tree $T$ is isometric to a real line then the associated $\Gamma$-limit group is abelian.  We have already covered this case, so we may assume $T$ is not isometric to a real line.  Therefore, by \ref{LinfProps}.(\ref{tripod}), \ref{LinfProps}.(\ref{SeqStable}) and \ref{LinfProps}.(\ref{fiStable}), $\{ f_i \}$ is a stable sequence and $K_\infty = \SK (f_i)$.  Hence the limit group $L_\infty$ is an algebraic $\Gamma$-limit group.

Conversely, suppose that $\{ h_i \co G \to \Gamma \}$ is a stable sequence of homomorphisms.  If the associated sequence of stretching factors $\{ \mu_j \}$ contains a bounded subsequence, then there is a subsequence $\{ h_{i'} \}$ of $\{ h_i \}$ so that $h_{i_1'}(G) \cong h_{i_2'}(G)$ for all $i_1'  , i_2' \in \{ h_{i'} \}$.  In this case, since $\{ h_i \}$ is a stable sequence, the associated algebraic $\Gamma$-limit group is isomorphic to a finitely generated subgroup of $\Gamma$.

Thus suppose that there is not a bounded subsequence of the $\{ \mu _j \}$.  In this case we can construct a limiting spaces $X_\omega$ and $\mathcal C_\infty$ and the associated $\R$-tree $T$.  If $T$ is isometric to a real line, we are done.  Otherwise since passing to a subsequence of a stable sequence does not change the stable kernel, we see again that $K_\infty = \SK (h_i)$, so that the algebraic $\Gamma$-limit group is a $\Gamma$-limit group.
\end{proof}

We now recall a topology on the set of finitely generated groups from \cite{CG} (see also \cite{Grig, Champ}).

\begin{definition}
A {\em marked group} $(G,\mathcal A)$ consists of a finitely generated group $G$ with an ordered generating set $\mathcal A = (a_1, \ldots , a_n )$.  Two marked groups $(G,\mathcal A)$ and $(G',\mathcal A')$ are {\em isomorphic} if the bijection taking $a_i$ to $a_i'$ for each $i$ induces an isomorphism between $G$ and $G'$.

For a fixed $n$, the set $\Gn$ consists of those marked groups $(G,\mathcal A)$ where $|\mathcal A| = n$.
\end{definition}

We now introduce a metric on $\Gn$.  First, we introduce the following abuse of notation.

\begin{convention} \label{Abuse}
Marked groups are always considered up to isomorphism of marked groups.  Thus for any marked groups $(G,\mathcal A)$ and $(G',\mathcal A')$ in $\Gn$, we identify an $\mathcal A$-word with the corresponding $\mathcal A'$-word under the canonical bijection induced by $a_i \to a_i'$, $i = 1 ,\ldots , n$.
\end{convention}

\begin{definition}
A {\em relation} in a marked group $(G,\mathcal A)$ is an $\mathcal A$-word representing the identity in $G$. Two marked groups $(G,\mathcal A)$ and $(G',\mathcal A')$ in $\Gn$ are at distance $e^{-d}$ from each other if they have the exactly the same relations of length at most $d$, but there is a relation of length $d+1$ which holds in one marked group but not the other.
\end{definition}

The following result is implicit in \cite{CG}.

\begin{proposition}
Let $G$ be a finitely generated group and $\Xi$ a finitely presented group.  Suppose that $\{ h_i \co G \to \Xi \}$ is a stable sequence, and $\{ g_1, \ldots , g_k \}$ is a generating set for $G$.  Then the marked group 
\[	\left(G/\SK \{ h_i \}, \left\{ g_1 \SK \{ h_i \}, \ldots , g_k \SK  \{ h_i \} \right\} \right)	,	\]
is a limit of marked groups $(G_i,\mathcal A_i)$ where each $G_i$ is a finitely generated subgroup of $\Xi$.

Conversely, if the marked group $(G,\mathcal A)$ is a limit of finitely generated subgroups of $\Xi$, then $(G,\mathcal A)$ is an algebraic $\Xi$-limit group.
\end{proposition}
\begin{proof}
Suppose that $\{ h_i \co G \to \Xi \}$ is a stable sequence, and $\{ g_1, \ldots , g_k \}$ is a generating set for $G$.  Consider the marked group
\[	\left(G/\SK \{ h_i \}, \left\{ g_1 \SK \{ h_i \}, \ldots , g_k \SK  \{ h_i \} \right\} \right)	.	\]
For each $n$, let $H_n = \langle h_n (g_1), \ldots , h_n(g_k) \rangle \leq \Xi$.  We consider the marked groups $\left( H_n , \left\{ h_n(g_1), \ldots , h_n(g_k) \right\} \right)$, and prove that they converge to $G/ \SK \{ h_i \}$ with the above marking.

Let $j \ge 1$ be arbitrary and let $W_j$ be the set of all words of length at most $j$ in the alphabet $\{ x_1^{\pm 1}, \ldots , x_k^{\pm 1} \}$.  Following Convention \ref{Abuse}, we interpret $W_j$ as words in the various generating sets without changing notation.  The set $W_j$ admits a decomposition into $T_j \cup N_j$, where $T_j$ are the words of length at most $j$ which are in $\SK \{ h_i \}$ and $N_j$ are the remaining words of length at most $j$.  Since $\{ h_i \}$ is a stable sequence, for each element $w \in T_j$, the element $h_n(w)$ is trivial for all but finitely many $n$, and for each $w \in N_j$ the element $h_n(w)$ is nontrivial for all but finitely many $n$.  

Thus, for all but finitely many $n$, the relations in $H_n$ of length at most $j$ are exactly the same as the relations of length at most $j$ in $G / \SK \{ h_i \}$.  Thus for all but finitely many $n$, the group $H_n$ with the given marking is at distance at most $e^{-j}$ from $G / \SK \{ h_i \}$ with the given marking.  This implies that the sequence $\{ H_n \}$ (with markings) converges to $G / \SK \{ h_i \}$ (with marking).

For the converse, suppose that $(G,\mathcal A)$ is a limit of a (convergent) sequence of marked finitely generated subgroups of $\Xi$.  Denote these subgroups by $(H_i,\mathcal A_i)$.  Note that $|\mathcal A_i|$ is fixed. Let $\mathcal A = \{ a_1 , \ldots , a_k \}$, let $\mathcal A_i = \{ b_{i,1}, \ldots , b_{i,k} \}$ and let $F$ be the free group on the set $\mathcal A$.  Define homomorphisms $h_i \co F \to \Xi$ by $h_i(a_j) = b_{i,j}$.  It is not difficult to see that $G \cong F / \SK \{ h_i \}$.
\end{proof}

For a group $H$, let $T_\forall(H)$ be the {\em universal theory} of $H$ -- the set of all universal sentences which are true in $H$ (see \cite{CG} or \cite{Paulin2} for the definition, we are interested only in its consequences).  The results in \cite{CG} now imply 

\begin{corollary}
Let $\Xi$ be a finitely presented group and suppose that $L$ is an algebraic $\Xi$-limit group.  Then $T_\forall(\Xi) \subseteq T_\forall (L)$.
\end{corollary}

\begin{corollary}
Suppose that $\Gamma$ is a torsion-free toral CAT$(0)$ group with isolated flats and that $L$ is a $\Gamma$-limit group.  Then
\begin{enumerate}
\item any finitely generated subgroup of $L$ is a $\Gamma$-limit group;
\item $L$ is torsion-free;
\item $L$ is commutative transitive; and
\item $L$ is CSA.
\end{enumerate}
\end{corollary}

Lemma \ref{SolAb} now implies the following:

\begin{corollary} \label{AbSubgps}
Let $\Gamma$ be a torsion-free toral CAT$(0)$ group with isolated flats and let $L$ be a  $\Gamma$-limit group. Every solvable subgroup of $L$ is abelian.
\end{corollary}

\subsection{The Main Theorem and conclusions}

Finally, we have:

\begin{theorem} \label{SplittingTheorem}
Suppose that $\Gamma$ is a torsion-free toral CAT$(0)$ group with isolated flats such that $\text{Out}(\Gamma)$ is infinite.  Then $\Gamma$ admits a nontrivial splitting over a finitely generated free abelian group.
\end{theorem}
\begin{proof}
Suppose that $\{ \phi_i \}$ is an infinite set of automorphisms of $\Gamma$ which belong to distinct conjugacy classes in $\text{Out}(\Gamma)$.  Then the construction from Sections \ref{Limits} and \ref{TreeSection} allows us to find an isometric action of $\Gamma$ on an $\R$-tree $T$ without global fixed points.  Pass to the subsequence $\{ f_i \}$ of $\{ \phi_i \}$ as in Lemma \ref{GromovTop}.

Suppose first that $T$ is isometric to a real line.  Then by Theorem \ref{LinfProps}.(\ref{Tline}) the group $f_i(\Gamma)$ is free abelian for all but finitely many $i$.  But $f_i(\Gamma) = \Gamma$, so $\Gamma$ must be free abelian in this case.  The theorem certainly holds for finitely generated free abelian groups.

Therefore, we may suppose that $T$ is not isometric to a real line.  In this case, since $K_\infty = \SK$ is trivial, the $\Gamma$-limit group $L_\infty$ is $\Gamma$ itself.  Then by Theorem 9.5 of \cite{BF}, the group $\Gamma$ splits over a group of the form $E$-by-cyclic, where $E$ fixes a non-degenerate segment of $T$.  The stabilisers in $\Gamma$ of non-degenerate segments are free abelian, by Lemma \ref{LinfProps}.(\ref{Stable}).  Hence the group of the form $E$-by-cyclic is soluble, and hence free abelian by  Lemma \ref{SolAb}.  Note that free abelian subgroups of $\Gamma$ are finitely generated.  This finishes the proof of the theorem.
\end{proof}

Suppose that $\Gamma$ is a CAT$(0)$ group and that $\text{Out}(\Gamma)$ is infinite.  Swarup asked (see \cite{BestvinaQuestions}, Q2.1) whether $\Gamma$ necessarily admits a Dehn twist of infinite order.  The above result shows that this is the case for the class of torsion-free toral CAT$(0)$ groups with isolated flats.  Swarup also asked whether there is an analog of the theorem of Rips and Sela that $\text{Out}(\Gamma)$ is virtually generated by Dehn twists.

In the subsequent work \cite{CWIFShort}, we will prove if $\Gamma$ is a torsion-free toral CAT$(0)$ group with isolated flats then $\text{Out}(\Gamma)$ is virtually generated by generalised Dehn twists (which take into account the existence of noncyclic abelian groups). 

It seems that the techniques developed here will be of little help in answering Swarup's question in the case of a general CAT$(0)$ group.

It is also clear from the above construction that if $L$ is a non-abelian freely indecomposable strict $\Gamma$-limit group then $L$ splits over an abelian group.  However, there is no reason to conclude that the edge group in this splitting is finitely generated.

It is straightforward to construct the canonical abelian JSJ decomposition of a strict $\Gamma$-limit group $L_\infty$, using acylindrical accessibility \cite{SelaAcyl}.\footnote{Note that strict $\Gamma$-limit groups need not be finitely presented, so the results of \cite{RipsSela} do not apply.}  However, for example, we do not yet know that the edge groups in the abelian JSJ decomposition of $L_\infty$ are finitely generated.  To prove that this is the case if $L_\infty$ is freely indecomposable and nonabelian involves the shortening argument of Sela, which we present for torsion-free toral CAT$(0)$ groups with isolated flats in future work.  The shortening argument allows us to prove that torsion-free toral CAT$(0)$ groups with isolated flats are Hopfian, and to construct Makanin-Razborov diagrams for these groups, thus partially answering a question of Sela (see \cite[I.8.(i), (ii), (iii)]{SelaProblems}).  This work will be undertaken in \cite{CWIFShort}.

\Addresses\recd

\end{document}

%% file: agtout.tex
\def\ifplaintex{\expandafter\ifx\csname documentclass\endcsname\relax}

\def\gtp{{\mathsurround=0pt\it $\cal G\mskip-2mu$eometry \&\ 
$\cal T\!\!$opology $\cal P\!$ublications}}  

\def\recd{{\small Received:\qua\receiveddate\ifx\reviseddate\relax
\else\qquad Revised:\qua\reviseddate\fi\par}} 

\def\lognumber#1{\def\thelognumber{#1}}
\def\volumenumber#1{\def\thevolumenumber{#1}}
\def\volumeyear#1{\def\thevolumeyear{#1}}
\def\papernumber#1{\def\thepapernumber{#1}}
\def\pagenumbers#1#2{\def\startpage{#1}\def\finishpage{#2}}
\def\published#1{\def\publishdate{#1}}

\def\received#1{\def\receiveddate{#1}}

\def\accepted#1{\def\accepteddate{#1}}

\long\def\asciiabstract#1{\long\def\theasciiabstract{#1}}
\def\asciikeywords#1{\def\theasciikeywords{#1}}

\let\\\par\let\thelognumber\relax\let\thevolumenumber\relax
\let\thepapernumber\relax\let\thevolumeyear\relax\let\startpage\relax
\let\finishpage\relax\let\publishdate\relax\let\receiveddate\relax
\let\reviseddate\relax\let\accepteddate\relax\let\theasciititle\relax
\let\theasciiauthors\relax
\let\theasciiabstract\relax\let\theasciikeywords\relax

\let\theasciiemail\relax

\ifplaintex
\font\logobig=cmssbx10 scaled 3836
\font\logomed=cmssbx10 scaled 2557
\else
\font\logobig=cmssbx10 scaled 4200
\font\logomed=cmssbx10 scaled 2800
\fi

\long\def\makeagttitle{   
\count0=\startpage
\agt\hfill      
\hbox to 45truept{\vbox to 0pt{\vglue -13truept{\logomed A\kern -.37em{\logobig 
T}\kern -.38em G}\vss}\hss}
\break
{\small Volume \thevolumenumber\ (\thevolumeyear)
\startpage--\finishpage\nl
Published: \publishdate}

\vglue .25truein

{\parskip=0pt\leftskip 0pt plus
1fil\def\\{\par\smallskip}{\Large\bf\thetitle}\par\medskip} \vglue
0.05truein

{\parskip=0pt\leftskip 0pt plus 1fil\def\\{\par}{\sc\theauthors}
\par\medskip}%
 
\vglue 0.03truein 

{\small\leftskip 25truept\rightskip 25truept{\bf Abstract}\stdspace\theabstract

{\bf AMS Classification}\stdspace\theprimaryclass
\ifx\thesecondaryclass\relax\else; \thesecondaryclass\fi\par
{\bf Keywords}\stdspace \thekeywords\par}\vglue 7truept

}   

\ifplaintex
\font\phead=cmsl9 scaled 950
\font\pnum=cmbx10 scaled 913
\font\pfoot=cmsl9 scaled 950
\headline{\vbox to 0pt{\vskip -4.5mm\line{\small\phead\ifnum
\count0=\startpage ISSN 1472-2739 (on-line) 1472-2747 (printed)
\hfill {\pnum\folio}\else\ifodd\count0\def\\{ }%
\ifx\theshorttitle\relax\thetitle\else\theshorttitle\fi\hfill{\pnum\folio}
\else\def\\{ and }{\pnum\folio}\hfill\ifx\theshortauthors\relax\theauthors
\else\theshortauthors\fi\fi\fi}\vss}}
\footline{\vbox to 0pt{\vglue 0mm\line{\small\pfoot\ifnum\count0=\startpage
\copyright\ \gtp\hfill\else
\agt, Volume \thevolumenumber\ (\thevolumeyear)\hfill\fi}\vss}}
\else
\font=cmsl9 scaled 1050
\font\lnum=cmbx10 
\font=cmsl9 scaled 1050
\makeatletter
\def\@oddhead{{\small\lhead\ifnum\count0=\startpage ISSN 1472-2739 
(on-line) 1472-2747 (printed)\hfill {\lnum\number\count0}\else\ifodd\count0
\def\\{ }\ifx\theshorttitle\relax \thetitle \else\theshorttitle\fi\hfill
{\lnum\number\count0}\else\def\\{ and }{\lnum\number\count0}
\hfill\ifx\theshortauthors\relax 
\theauthors\else\theshortauthors\fi\fi\fi}}\def\@evenhead{\@oddhead}
\def\@oddfoot{\small\lfoot\ifnum\count0=\startpage\copyright\ \gtp\hfill\else
\agt, Volume \thevolumenumber\ (\thevolumeyear)\hfill\fi}
\def\@evenfoot{\@oddfoot}
\makeatother
\fi
\let\maketitlepage\makeagttitle
\let\makeshorttitle\maketitlepage
\let\maketitle\maketitlepage

\newwrite\gtoutfile
\long\gdef\makeheadfile{  
{\def\\{, }\def\s{ }
\immediate\openout\gtoutfile head.xxx
\immediate\write\gtoutfile{Proxy-for: \ifx\theasciiauthors\relax
\theauthors\else\theasciiauthors\fi\s<\ifx\theasciiemail\relax\theemail\else\theasciiemail\fi>}
\immediate\write\gtoutfile{\noexpand\\}
\immediate\write\gtoutfile{Authors: \ifx\theasciiauthors\relax
\theauthors\else\theasciiauthors\fi}
{\def\\{ }\immediate\write\gtoutfile{Title: \ifx\theasciititle\relax
\thetitle\else\theasciititle\fi}}
\immediate\write\gtoutfile{Subj-class: GT or SG, GR etc}
\immediate\write\gtoutfile{MSC-class: \theprimaryclass\ifx\thesecondaryclass\relax\else, \thesecondaryclass\fi}
\immediate\write\gtoutfile{Journal-ref: Algebr. Geom. Topol. \thevolumenumber\s
(\thevolumeyear) \startpage-\finishpage}
\immediate\write\gtoutfile{Comments: Published by Algebraic and
Geometric Topology at}
\immediate\write\gtoutfile{\s\s\s  http://www.maths.warwick.ac.uk/agt/AGTVol\thevolumenumber/agt-\thevolumenumber-\thepapernumber.abs.html}
\immediate\write\gtoutfile{\noexpand\\}
\immediate\write\gtoutfile{}
\ifx\theasciiabstract\relax
\immediate\write\gtoutfile{\theabstract}\else
\immediate\write\gtoutfile{\theasciiabstract}\fi
\immediate\write\gtoutfile{}
\immediate\write\gtoutfile{\noexpand\\}
\immediate\write\gtoutfile{}
\immediate\closeout\gtoutfile}}  

\def\maketitlepage{\makeagttitle\makeheadfile}
\let\makeshorttitle\maketitlepage
\let\maketitle\maketitlepage

%% file: RelThin.eps_t
\begin{picture}(0,0)%
\includegraphics{RelThin.eps}%
\end{picture}%
\setlength{\unitlength}{3947sp}%
\begingroup\makeatletter\ifx\SetFigFont\undefined%
\gdef\SetFigFont#1#2#3#4#5{%
  \reset@font\fontsize{#1}{#2pt}%
  \fontfamily{#3}\fontseries{#4}\fontshape{#5}%
  \selectfont}%
\fi\endgroup%
\begin{picture}(2424,2615)(4414,-2923)
\put(5551,-1861){\makebox(0,0)[lb]{\smash{{\SetFigFont{12}{14.4}{\rmdefault}{\mddefault}{\updefault}{\color[rgb]{0,0,0}$E$}%
}}}}
\end{picture}%

%% file: Approx.eps_t
\begin{picture}(0,0)%
\includegraphics{Approx.eps}%
\end{picture}%
\setlength{\unitlength}{3947sp}%
\begingroup\makeatletter\ifx\SetFigFont\undefined%
\gdef\SetFigFont#1#2#3#4#5{%
  \reset@font\fontsize{#1}{#2pt}%
  \fontfamily{#3}\fontseries{#4}\fontshape{#5}%
  \selectfont}%
\fi\endgroup%
\begin{picture}(2823,2271)(2026,-2491)
\put(3190,-282){\makebox(0,0)[lb]{\smash{{\SetFigFont{6}{7.2}{\rmdefault}{\mddefault}{\updefault}{\color[rgb]{0,0,0}$\A_k$}%
}}}}
\put(3297,-1128){\makebox(0,0)[lb]{\smash{{\SetFigFont{6}{7.2}{\rmdefault}{\mddefault}{\updefault}{\color[rgb]{0,0,0}$\A_k'$}%
}}}}
\put(2026,-2150){\makebox(0,0)[lb]{\smash{{\SetFigFont{6}{7.2}{\rmdefault}{\mddefault}{\updefault}{\color[rgb]{0,0,0}$\B_k$}%
}}}}
\put(2802,-1833){\makebox(0,0)[lb]{\smash{{\SetFigFont{6}{7.2}{\rmdefault}{\mddefault}{\updefault}{\color[rgb]{0,0,0}$\B_k'$}%
}}}}
\put(4849,-2468){\makebox(0,0)[lb]{\smash{{\SetFigFont{6}{7.2}{\rmdefault}{\mddefault}{\updefault}{\color[rgb]{0,0,0}$\C_k$}%
}}}}
\put(3650,-1657){\makebox(0,0)[lb]{\smash{{\SetFigFont{6}{7.2}{\rmdefault}{\mddefault}{\updefault}{\color[rgb]{0,0,0}$\C_k'$}%
}}}}
\put(3155,-1552){\makebox(0,0)[lb]{\smash{{\SetFigFont{6}{7.2}{\rmdefault}{\mddefault}{\updefault}{\color[rgb]{0,0,0}$\overline{N}_k$}%
}}}}
\put(2132,-1198){\makebox(0,0)[lb]{\smash{{\SetFigFont{6}{7.2}{\rmdefault}{\mddefault}{\updefault}{\color[rgb]{0,0,0}$E$}%
}}}}
\end{picture}%